\documentclass{amsart}
\usepackage{latexsym,amsxtra,amscd,ifthen}
\usepackage{amsfonts}
\usepackage{verbatim}
\usepackage{amsmath}
\usepackage{amsthm}
\usepackage{amssymb}

\numberwithin{equation}{section}

\theoremstyle{plain}
\newtheorem{theorem}{Theorem}[section]
\newtheorem{lemma}[theorem]{Lemma}
\newtheorem{proposition}[theorem]{Proposition}
\newtheorem{corollary}[theorem]{Corollary}
\newtheorem{conjecture}[theorem]{Conjecture}

\theoremstyle{definition}
\newtheorem{definition}[theorem]{Definition}
\newtheorem{example}[theorem]{Example}

\newtheorem{remark}[theorem]{Remark}

\makeatletter              
\let\c@equation\c@theorem  
\makeatother

\newcommand{\wt}{\widetilde}
\newcommand{\lra}{\longrightarrow}
\newcommand{\Alb}{\operatorname{Alb}}

\newcommand{\Num}{\operatorname{Num}}
\newcommand{\pt}{\operatorname{pt}}
\DeclareMathOperator{\cdim}{cd}
\DeclareMathOperator{\caphom}{Hom}
\DeclareMathOperator{\GL}{GL}
\DeclareMathOperator{\GK}{GK}
\DeclareMathOperator{\uExt}{\underline{Ext}}
\DeclareMathOperator{\uHom}{\underline{Hom}}
\DeclareMathOperator{\uotimes}{\underline{\otimes}}
\DeclareMathOperator{\Ext}{Ext}

\DeclareMathOperator{\Pic}{Pic}
\DeclareMathOperator{\HB}{H}
\DeclareMathOperator{\Bir}{Bir}
\DeclareMathOperator{\HS}{h}
\DeclareMathOperator{\ann}{ann}

\DeclareMathOperator{\rann}{r.ann}
\DeclareMathOperator{\cd}{cd}
\DeclareMathOperator{\QGr}{QGr}
\DeclareMathOperator{\Qgr}{QGr}
\DeclareMathOperator{\gr}{gr}
\DeclareMathOperator{\Gr}{Gr}
\DeclareMathOperator{\qgr}{qgr}

\DeclareMathOperator{\proj}{proj}
\DeclareMathOperator{\tors}{tors}
\DeclareMathOperator{\Tors}{Tors}

\DeclareMathOperator{\Proj}{Proj}

\DeclareMathOperator{\Aut}{Aut}
\DeclareMathOperator{\Kdim}{Kdim}
\DeclareMathOperator{\GKdim}{GKdim}
\DeclareMathOperator{\End}{End}
\DeclareMathOperator{\R}{R}
\DeclareMathOperator{\D}{D}
\DeclareMathOperator{\Hom}{Hom}
\DeclareMathOperator{\RHom}{RHom}

\DeclareMathOperator{\im}{im}

\newcommand{\fm}{\mathfrak{m}}

\newcommand{\cal}{\mathcal}
\newcommand{\mf}{\mathfrak}
\newcommand{\mc}{\mathcal}
\newcommand{\mb}{\mathbb}
\newcommand{\id}{\operatorname{id}}

\newcommand{\ch}{\operatorname{char}}
\newcommand{\Qch}{\operatorname{Qch}}

\begin{document}

\title{Projectively Simple Rings}

\author{Z. Reichstein, D. Rogalski and J. J. Zhang}

\address{(Reichstein) Department of Mathematics, 
University of British Columbia, 
Vancouver, British Columbia, Canada V6T 1Z2}
\thanks{Z. Reichstein was partially supported by an NSERC research grant}

\email{reichst@math.ubc.ca}

\address{(Rogalski) Department of Mathematics, MIT,
Cambridge, MA 02139-4307, USA}
\thanks{D. Rogalski was partially supported by NSF grant DMS-0202479}

\email{rogalski@math.mit.edu}

\address{(Zhang) Department of Mathematics, Box 354350,
University of Washington, Seattle, WA 98195, USA}
\thanks{J. J. Zhang was partially supported by NSF grant DMS-0245420}

\email{zhang@math.washington.edu}

\begin{abstract}
An infinite-dimensional $\mb{N}$-graded $k$-algebra $A$ is called 
\emph{projectively simple} if $\dim_k A/I < \infty$
for every nonzero two-sided ideal $I \subset A$. 
We show that if a projectively simple
ring $A$ is strongly noetherian, is generated in degree $1$, and has a
point module,
then $A$ is equal in large degree to a twisted homogeneous coordinate ring
$B = B(X,
\mc{L}, \sigma)$.  Here $X$ is a smooth projective variety, $\sigma$ is an
automorphism of $X$ with no proper $\sigma$-invariant subvariety (we call
such
automorphisms \emph{wild}), and $\mc{L}$ is a $\sigma$-ample line bundle.
We
conjecture that if $X$ admits a wild automorphism then every irreducible
component of
$X$ is an abelian variety.  We prove several results in support of this
conjecture;
in particular, we show that the conjecture is true if $\dim X \leq 2$.  In
the case
where $X$ is an abelian variety, we describe all wild automorphisms of
$X$. Finally,
we show that if $A$ is projectively simple and admits a balanced dualizing
complex,
then $\proj A$ is Cohen-Macaulay and Gorenstein.
\end{abstract}

\subjclass[2000]{16W50, 14A22, 14J50, 14K05}


\keywords{Graded ring, simple ring, just infinite dimensional algebra,
twisted homogeneous coordinate ring,
projective variety, algebraic surface, wild automorphism, abelian
variety, dualizing complex}

\maketitle

\tableofcontents

\setcounter{section}{-1}
\section{Introduction}
\label{sec0}

Let $k$ be a field.
The present study grew out of the following question:  What are the simple
graded $k$-algebras $A = \bigoplus_{n = 0}^{\infty} A_n$?
Technically, such an algebra cannot be simple, since it always has
ideals of the form $A_{\geq n}$ for each $n \geq 0$. Thus we are
led to the following natural definition:
we call $A$ \emph{projectively simple}
if $\dim_k A = \infty$ and every nonzero graded ideal
of $A$ has finite codimension in $A$. If the reference to $k$ is
clear from the context, we will sometimes refer to $A$ as a
projectively simple graded ring.

Although it is easy to see that graded prime algebras of GK-dimension $1$
which satisfy a polynomial identity are projectively simple, it
is not immediately obvious that there are any examples
more interesting than these.  In fact, one interesting example
has already appeared prominently in the literature: If $S$ is the Sklyanin
algebra of
dimension $3$ (as defined in \cite[Example 7.3]{SV}), then $S$ has a
central element
$g$ such that the factor ring $B = S/gS$ is a projectively simple domain
of
GK-dimension $2$.  In this case $B$ may be constructed as a twisted
homogeneous
coordinate ring $B = B(E, \mc{L}, \sigma)$, where $E$ is an elliptic curve
and
$\sigma$ is a translation automorphism of $E$.  Thus it is natural to look
for other
examples of projectively simple rings which also arise as twisted
homogeneous
coordinate rings; let us briefly review this construction.

Given $X$ a projective scheme, $\mc{L}$ a line bundle on $X$, and $\sigma$
an
automorphism of $X$, set 
$\mc{L}_0=\mc{O}_X$ and 
$\mc{L}_n = \mc{L} \otimes \sigma^{*} \mc{L}
\otimes \dots
\otimes (\sigma^{n-1})^{*}\mc{L}$ for each $n \geq 1$. 
The \emph{twisted
homogeneous
coordinate ring} associated to this data is the vector space $B(X, \mc{L},
\sigma) =
\bigoplus_{n = 0}^{\infty} \HB^0(\mc{L}_n)$, which has a natural graded
ring
structure. As is true even in the commutative case where $\sigma = Id$,
the rings
obtained this way are typically well behaved only with additional
assumptions on the
sheaf $\mc{L}$. In particular, if $\mc{L}$ is \emph{$\sigma$-ample} (see
\S\ref{xxsec2} below), then $B = B(X, \mc{L}, \sigma)$ is noetherian and
its
properties are closely related to the geometry of $X$.  In this case one
can describe
the ideals of $B$ geometrically and thus prove the following:
\begin{proposition} (Proposition~\ref{prop2.2})
\label{prop0.1}
Let $X$ be a projective scheme, $\sigma$ an automorphism of $X$, and
$\mc{L}$ a $\sigma$-ample line bundle.
Then $B = B(X, \mc{L}, \sigma)$ is projectively simple if and only if
\[ \sigma(Y)\neq Y\ \text{for every nonempty reduced closed subscheme}\
Y\subsetneq X.
\]
\end{proposition}

In the sequel we shall refer to $\sigma$ satisfying this condition as a
{\em wild} automorphism of $X$.
If we seek explicit examples of projectively simple twisted homogeneous
coordinate
rings, then Proposition~\ref{prop0.1} inevitably leads us to the following
purely
geometric questions:

\smallskip
(i) Can we classify all pairs $(X, \sigma)$, where
$X$ is a projective variety and $\sigma$ is a wild automorphism of $X$?

\smallskip
(ii) For each $(X, \sigma)$ as above, can we find $\sigma$-ample
line bundles on $X$?

\smallskip
In the case where $X$ is an abelian variety, we
obtain the following complete answer to these questions.
Recall that every automorphism $\sigma$ of an abelian variety
$X$ can be written in the form $\sigma = T_b \cdot \alpha$, where
$T_b: x \mapsto x+ b$ is a translation by some $b \in X$ and $\alpha$
is an automorphism of $X$ preserving the group structure; see, for
example, \cite[Theorem 4, p. 24]{La1}.

\begin{theorem} (Theorem~\ref{thm6.4}, Theorem~\ref{prop7.4})
\label{thm0.2}
Let $X$ be an abelian variety over
an algebraically closed field $k$ of characteristic zero.
\begin{enumerate}
\item Suppose that $T_b: X \lra X$ is the translation
automorphism by $b \in X$, and $\alpha$ is an automorphism
of $X$ preserving the group structure.
Then $\sigma = T_b \cdot \alpha$
is a wild automorphism of $X$ if and only if $\beta = \alpha - Id$ is
nilpotent and $b$ generates $X/\beta(X)$.
\item If $\sigma$ is a wild automorphism of $X$, then
any ample invertible sheaf on $X$ is $\sigma$-ample.
\end{enumerate}
\end{theorem}

Given an abelian variety $X$, we shall see in \S\ref{sec6} 
below that it is easy to find many automorphisms $\sigma$ satisfying
the conditions in Theorem~\ref{thm0.2}(a).  Then
as $\mc{L}$ varies over all ample sheaves on $X$, by
Theorem~\ref{thm0.2}(b) and Proposition~\ref{prop0.1}
one gets many projectively simple twisted homogeneous coordinate rings
$B(X, \mc{L}, \sigma)$. When $k$ is uncountable, for every integer 
$n>0$, there are  many noetherian projectively simple rings of 
GK-dimension $n$.

We conjecture that there are no examples of wild automorphisms,
other than those described in Theorem~\ref{thm0.2}.
More precisely,

\begin{conjecture}
\label{con0.3}
If an irreducible projective variety $X$ admits a wild automorphism
then $X$ is an abelian variety.
\end{conjecture}

In Sections~\ref{xxsec3} -- \ref{sec5d} we prove a number of results in
support in
Conjecture~\ref{con0.3}. In particular, 
we show that this conjecture is true if $\dim X \le 2$; see
Theorem~\ref{thm5.12}. 

Having constructed some prototypical examples, we would like
to say more about the structure of general projectively
simple rings.  Because twisted homogeneous coordinate rings
are quite special among all graded rings, one might expect
that the examples we have constructed so far are also very special.
Rather surprisingly, the following theorem states that under certain
hypotheses which are natural in the theory of noncommutative
projective geometry, these are essentially the only possible examples.

\begin{theorem} (Theorem~\ref{thm2.5})
\label{thm0.4}
Let $k$ be an algebraically closed field and let $A$ be a
projectively simple noetherian $k$-algebra. Suppose that $A$
is strongly noetherian, generated in degree 1 and has a point module.
Then there is an injective homomorphism $A \hookrightarrow B$ 
of graded algebras, such that $\dim_k \,  B/A < \infty$ 
and $B = B(X,{\mathcal L},\sigma)$ is a projectively simple
twisted homogeneous coordinate ring for some smooth projective 
variety $X$ with a wild automorphism $\sigma$ and
a $\sigma$-ample line bundle ${\mathcal L}$. 
\end{theorem}
It would be very interesting to know what kinds of projectively simple
rings can
appear if the various hypotheses of Theorem~\ref{thm0.4} are relaxed.  In
Example~\ref{ex2.6} we show that there do exist projectively simple
rings which are not strongly noetherian.

Our final class of results concerns the noncommutative projective
scheme $\Proj A$ associated to a projectively simple ring $A$, and
takes its inspiration from some known results about ungraded simple
rings. It is shown in \cite{YZ2} that a noetherian simple ring $A$ has
finite injective dimension provided that $A$ admits a dualizing complex.
This last condition is a natural one which holds for many important
classes of rings such as the universal enveloping algebras of finite
dimensional Lie algebras and factor rings of Artin-Schelter regular rings.
Following the ideas in~\cite{YZ2}, we show that some similar homological
results are true in the graded setting.

\begin{theorem} (Proposition~\ref{prop2.1}, Theorem~\ref{thm3.6})
\label{thm0.5}
Let $A$ be a projectively simple noetherian connected
graded algebra with a balanced dualizing complex. Then
\begin{enumerate}
\item $\Proj A$ is classically Cohen-Macaulay, namely, the dualizing
complex
for $\Proj A$ is isomorphic to $\omega[n]$ for a graded $A$-bimodule
object
$\omega$ in $\Proj A$.
\item $\Proj A$ is Gorenstein, namely, the dualizing bimodule $\omega$ in
part
(a) is invertible. As a consequence, the structure sheaf ${\mathcal A}$
has
finite injective dimension.
\end{enumerate}
\end{theorem}

If $A$ satisfies the hypotheses of Theorem~\ref{thm0.5} and also $A =
B(X,{\mathcal L},\sigma)$ for some $\sigma$-ample $\mc{L}$, then the
conclusions of
Theorem~\ref{thm0.5} are immediate, since $\Qgr A$ is equivalent to the
category of quasi-coherent sheaves on $X$ and $X$ is smooth by
Lemma~\ref{lem2.3}(b) below.  So Theorem~\ref{thm0.5} is aimed 
primarily at projectively simple rings that
do not satisfy the conditions of Theorem~\ref{thm0.4}.

To conclude this introduction, we would like to say a bit about our
motivation for
introducing and studying projectively simple rings. Let us describe a
potential application to the theory of GK-dimension; see \cite{KL} for an
introduction to this subject. It is known that no algebra may have a
non-integer value of
GK-dimension between $0$ and $2$.  Moreover, it is conjectured that there
does not
exist a noetherian connected graded domain with GK-dimension strictly
between $2$
and $3$ \cite[p. 2]{AS1}. Note that if $A$ is a Goldie prime ring and $I$
is a
nonzero ideal of $A$, then $\GK(A/I) \leq \GK(A) - 1$ \cite[Proposition
3.15]{KL}.  Thus
if $A$ is a connected graded noetherian prime ring with $2 < \GK(A) < 3$,
then
either $A$ has a height one prime $P$ with $\GK(A/P) = 1$, or else $A$ is
projectively
simple.  A proof that projectively simple rings have integer GK-dimension
would
eliminate one of these cases, thus making progress towards the proof of
the full conjecture.\footnote{A proof of this conjecture was recently announced by 
A. Smoktunowicz.}
The results of this paper do at least show that a
projectively simple ring
which satisfies the hypotheses of Theorem~\ref{thm0.4} must have integer
GK-dimension, since this is true of twisted homogeneous coordinate rings.

Another application is to the classification of graded rings of low
GK-dimension. Artin and Stafford have classified semiprime graded rings
of GK-dimension $2$ in terms of geometric data in \cite{AS1}, \cite{AS2}.
The classification of rings of GK-dimension $3$, which
correspond to noncommutative surfaces, is a subject of much current
interest.  Those rings of GK-dimension $3$ which are also
projectively simple represent a special subclass.
If such a ring satisfies all of the hypotheses of
Theorem~\ref{thm0.4} above, then it must be
a twisted homogeneous coordinate ring in large degree.
We classify all such twisted homogeneous coordinate rings in
Proposition~\ref{prop8.1}.

\section{Elementary properties of projectively simple rings}
\label{xxsec1}

Throughout this paper $k$ is a commutative base field, and all rings
will be $k$-algebras. In Section \ref{xxsec2} we will assume that $k$ is
algebraically closed, and is Sections \ref{xxsec3}-\ref{sec8} 
we will assume that
$k$ is algebraically closed and that $\ch k = 0$. An algebra $A$ is
{\it ${\mathbb N}$-graded} (or {\it graded}) if $A=\bigoplus_{i\geq 0}
A_i$ with $1\in A_0$ and $A_i A_j \subseteq A_{i+j}$ for all $i,j \geq 0$.
The graded algebra $A$ is {\it locally finite} if $\dim_k A_i<\infty$
for all $i$.  All algebras $A$ in this article will be graded and locally
finite, except when we consider localizations of $A$ and in other obvious
situations. If $A_0=k$, then $A$ is called {\it connected graded}.
Let $\fm$ denote the graded ideal $A_{\geq 1}$.

We recall from the introduction the property which is the subject of this
article:
\begin{definition}
\label{def1.1}
A locally finite graded algebra $A$ is called {\it projectively simple}
if $\dim_k A = \infty$ and $\dim_k A/I <\infty$ for every nonzero graded
ideal $I$ of $A$.
\end{definition}

\noindent
It is useful to notice that since $A$ is locally finite, the condition
$\dim_k A/I
<\infty$ for a graded ideal $I$ is equivalent to the condition $I\supset
A_{\geq n}$
for some $n$.

\begin{remark} 
An infinite group is called \emph{just infinite} if every
nontrivial normal subgroup has finite index. A $k$-algebra 
$A$ is called \emph{just infinite dimensional} if every nonzero
two-sided ideal has finite codimension in $A$
\cite{FaS, PT, Si1, Si2, Vi}.  Thus a projectively simple
algebra may be viewed as a graded counterpart of a just
infinite dimensional algebra. We are grateful to Lance
Small for bringing this notion and the above references
to our attention.
\end{remark}

Our first result summarizes some easy observations about projectively
simple rings. The proofs are straightforward, and we omit them.

\begin{lemma}
\label{lem1.2}
Let $A$ be a graded $k$-algebra such that $\fm\neq 0$.
\begin{enumerate}
\item
If $A$ is projectively simple, then $A$ is a finitely generated
$k$-algebra.
\item
If $A$ is projectively simple, then $A$ is prime. If $A$ is also
connected,
then the only nonzero graded prime ideal of $A$ is $\fm$.
\item
If $A$ is PI and finitely generated, then $A$ is projectively
simple if and only if $A$ is prime of GK-dimension 1.
\item
If $A$ is projectively simple, then $\dim_k A/J<\infty$ for every (not
necessarily graded) nonzero ideal $J$ of $A$. 
In other words, $A$ is a just infinite dimensional algebra.
\end{enumerate}
\end{lemma}

Before pursuing further the abstract notion of a projectively simple ring,
we should
note that there do exist examples besides the trivial ones in
Lemma~\ref{lem1.2}(c)
above. Let us mention here just one important example, which will turn out
to be the
model for all of our subsequent examples.

\begin{example}
\label{ex1.3}
Let $E$ be an elliptic curve and let $\sigma$ be an translation
automorphism of $E$ given by the rule $x \mapsto x+a$ where $a$ is a
point on $E$ of infinite order. Let ${\mathcal L}$ be a very ample
line bundle of $E$. Then the twisted homogeneous coordinate ring
$B(E, {\mathcal L} ,\sigma)$ (see \S\ref{xxsec2} for the definition)
is a projectively simple ring of GK-dimension 2 (see
\cite[6.5(ii)]{AS1} for a proof).
\end{example}

The following lemma is a graded analogue of \cite[Lemma 2.1]{FaS}, 
due to Small.

\begin{lemma}
\label{lem1.4}
Let $A$ be finitely generated, connected graded and infinite
dimensional over $k$.  Then there is a graded ideal $J\subset A$ such
that $A/J$ is projectively simple.
\end{lemma}

\begin{proof}
Let $\Phi$ be the set of all graded ideals $I\subset A$ such that
$\dim_k A/I=\infty$. Suppose that 
$\{I_\alpha\}_{\alpha\in S} \subset \Phi$ is an ascending chain
of (possibly
uncountably many) ideals. We claim that the union
$$U=\bigcup_{\alpha\in S} I_{\alpha}$$
is in $\Phi$. Suppose to the contrary that $U$ is not in $\Phi$. Then $U$
contains $A_{\geq n}$ for some $n$. Let $A$ be generated by elements
of degree no more than $d$ and let $V$ be the finite dimensional vector
space $\bigoplus_{i=0}^d A_{n+i}$. Then $V\subset U$ implies that
$V\subset I_{\alpha_0}$ for some $\alpha_0$ since
$\{I_\alpha\}_{\alpha\in S}$ is an ascending chain. Hence $I_{\alpha_0}$
contains $A_{\geq n}$ since $A$ is generated by elements of degree no
more than $d$. This yields a contradiction, whence $U$ is in $\Phi$. By
Zorn's lemma, the set $\Phi$ has a maximal object, say $J$. Then $A/J$
is projectively simple.
\end{proof}

The following result is a variant of a theorem of Farkas and Small 
\cite[Theorem 2.2]{FaS}.

\begin{proposition}
Let $A$ be a noetherian projectively simple algebra over an 
uncountable field $k$. If $\GKdim A>1$, then $A$ is primitive.
\end{proposition}

\begin{proof} By Lemma \ref{lem1.2}(b), $A$ is prime. By 
\cite[Theorem 0.4]{ASZ}, $A$ is Jacobson. Hence $A$ is
semiprimitive. Since $\GKdim A>1$, $A$ is not PI; see 
Lemma~\ref{lem1.2}(c). 
By \cite[Theorem 2.2]{FaS}, $A$ is primitive.
\end{proof}

One might ask if Lemma~\ref{lem1.2}(a) could be strengthened to say that a
projectively simple ring is necessarily noetherian. Using
Lemma~\ref{lem1.4} we can see that the answer is no. The following
construction is well-known.

\begin{example}
\label{ex1.6}
By the Golod-Shafarevitch construction, there is a 
finitely generated 
connected graded $k$-algebra $A$ such that $\dim_k A=\infty$ and 
every element in $A_{\geq 1}$ is nilpotent (see \cite[Chapter 8]{He}
when $k$ is countable and 
see \cite{FiS} when $k$ is an arbitrary field). 
By Lemma~\ref{lem1.4}, there is a graded ideal $J$
such that $B:=A/J$ is projectively simple, whence prime. Clearly
$B_{\geq 1}$ is a nil ideal, in other words every element is nilpotent,
but $B_{\geq 1}$ is
not nilpotent since $B$ is prime. By \cite[2.3.7]{MR}, $B$ is neither left
nor right noetherian.
\end{example}

Next we will show that projective simplicity is preserved under base
field extensions.

\begin{lemma}
\label{lem1.5}
Let $k \subset L$ be an extension of fields, with $k$ algebraically
closed.  Then $A$ is a projectively simple $k$-algebra if and only if
$A \otimes_k L$ is a projectively simple $L$-algebra.
\end{lemma}
\begin{proof}
If $R$ is any commutative $k$-algebra and $V$ is any vector space over
$k$, we write
$V_R$ for $V\otimes_k R$.  Suppose that $A$ is not projectively simple,
but rather
contains a nonzero ideal $I$ with $\dim_k A/I = \infty$. Then $I_L$ is a
nonzero
ideal of $A_L$ and $\dim_L A_L/I_L = \infty$, so $A_L$ is not projectively
simple.

Now assume that $A_L$ is not projectively simple, and let us show that $A$
is not
projectively simple.   By assumption $A_L$ has some nonzero ideal $I$ such
that
$A_L/I$ is not finite-dimensional. We may assume that $I=(g)$ is generated
by one
nonzero homogeneous element $g\in A_L$. Write $g=\sum_{i=0}^q b_i x_i$
where $0\neq
b_i\in L$ and $x_i\in A$. We may choose $q$ as small as possible so that
$\{x_0,\dots,x_q\}$ are linearly independent over $k$. Replacing $g$ by
$b_0^{-1}g$
we may assume that $b_0=1$. Let $R = k[b_1, \dots, b_q] \subset L$ be the
commutative
affine $k$-algebra generated by the $\{b_i \}$. Let $J$ be the ideal of
$A_R$
generated by $g$. Then $(A_R/J) \otimes_R L \cong A_L/I$.  Necessarily the
degree $n$
part $(A_R/J)_n$ is not $R$-torsion for all $n$ in an infinite subset
$\Phi\subset
\mathbb{Z}$.  Now let $\theta: R \to k$ be any $k$-algebra homomorphism
(some such
exists by the Nullstellensatz since $k$ is algebraically closed), and
extend it to a
map $\theta: A_R \to A_k = A$.  Then letting $k$ be an $R$-module via
$\theta$, we
have
\[
A' = (A_R/J) \otimes_R k \cong A/
(\theta(J)).
\]
For each $n\in \Phi$, $A'_n\neq 0$ since $(A_R/J)_n$ is not $R$-torsion.
Hence $A'$
has infinite $k$-dimension.  Now $\theta(g)=x_0+\sum_{i=1}^q
\theta(b_i)x_i$ is
nonzero in $\theta(J)$ since $\{x_0,\cdots,x_q\}$ is linearly independent
over $k$.
Therefore $A$ is not projectively simple.
\end{proof}

Now we discuss some further definitions which will be useful in the
remainder of this
section, and in \S\ref{xxsec10}-\ref{xxsec11}.
Let $A$ be a noetherian graded ring,
and let $M$ and $N$ be graded right $A$-modules.  The group of module
homomorphisms
preserving degree is denoted $\Hom_A(M,N)$. Let $N(n)$ denote the $n$th
degree shift
of the module $N$. We write
$$\uHom_A(M,N) = \bigoplus_{n \in \mb{Z}} \Hom_A(M, N(n)),$$
which is the group of all module homomorphisms from $M$ to $N$
in case $M$ is finitely generated.

We say that a graded right module $M$ is {\it torsion} if for every
$x\in M$ there is an $n$ such that $x A_{\geq n}=0$.  A noetherian
module is torsion if and only if it is finite dimensional over $k$.
We will also use a different kind of torsion property.  A right
$A$-module $M$ is called {\it Goldie torsion} if for every element
$x\in M$, the right annihilator $\rann(x)$ is an essential right
ideal of $A$.  Let $\Kdim M$ denote the Krull dimension of $M$;
if $A$ is prime then $\Kdim M < \Kdim A$ if and only if $M$ is Goldie
torsion. We call $M$ {\it torsionfree} (respectively, \emph{Goldie
torsionfree}) if it does not contain a nonzero torsion (respectively,
Goldie torsion) submodule.

Let $A$ and $B$ be graded rings.  An $(A, B)$-bimodule is called
{\it noetherian} if it is noetherian on both sides.  Note that if
$M$ is a noetherian $(A, B)$-bimodule, then the largest torsion left
submodule $\tau_A(M)$ of $M$ and the largest torsion right submodule
$\tau_B(M)$ of $M$ coincide, both being equal to the largest
finite-dimensional sub-bimodule of $M$. So we simply write $\tau(M)$
for this module and call it the torsion submodule of $M$. The bimodule
$M/\tau(M)$ is torsionfree on both sides.

The following useful result shows how the
projectively simple property gives strong information about the
structure of graded bimodules.
In particular, the two different notions of torsion defined above actually
coincide for noetherian bimodules over projectively simple rings.

\begin{lemma}
\label{lem1.7}
Let $A$ be a noetherian projectively simple ring.
\begin{enumerate}
\item
Let $B$ be any graded ring. If $M$ is a noetherian $(B,A)$-bimodule
such that $M_A$ is Goldie torsion, then $M$ is finite dimensional,
whence torsion over $A$.
\item
Let $M$ and $N$ be noetherian graded $(A, A)$-bimodules such that $M_A$
and $N_A$ are not torsion.  Then $\dim_k \uHom_A(M_A,N_A) = \infty$.
\end{enumerate}
\end{lemma}

\begin{proof}
(a) Since $_B M$ is noetherian, $M=\sum_i Bx_i$ for a finite set
of elements $\{x_i\}\subset M$. Thus
$$I:=\rann(M_A)=\bigcap_i \rann(x_i)\neq 0,$$
where the final inequality follows from the fact that every $\rann(x_i)$
is essential. Thus $M_A$ is finitely generated over a finite dimensional
algebra $A/I$. Therefore $M$ is finite dimensional.

(b) By part (a), $M_A$ and $N_A$ are not Goldie torsion. Since
$A$ is a prime Goldie ring, by replacing $M_A$ by
a factor module and $N_A$ by a submodule, we may assume that
both $M$ and $N$
are uniform Goldie torsionfree graded right $A$-modules, or homogeneous
right ideals of $A$.
Since $A$ is prime, $N_{\geq n}M\neq 0$ for all $n$. Hence for each $n$
there is a homogeneous element $x\in N_{\geq n}$  such that the left
multiplication
$$l_x: m\to xm\quad \text{for all } m\in M$$
is a nonzero element in $\uHom_A(M_A,N_A)$. Thus (b) follows.
\end{proof}

In the final result of this section, we note some circumstances under
which the
property of projective simplicity passes from one algebra to closely
related algebras.

\begin{lemma}
\label{lem1.8}
Let $A$ be a noetherian graded prime ring.
\begin{enumerate}
\item
Suppose that $B$ is a graded subring of $A$ such that $\dim_k A/B<\infty$.
Then $A$ is projectively simple if and only if $B$ is.
\item
Suppose that $B$ is a projectively simple graded subring of a graded
Goldie
prime ring $A$ such that $A_B$ is
finitely generated. Then $A$ is projectively simple.
\item
Suppose that $A$ is projectively simple, and let $B = A^{(n)}$ be the
$n$th
Veronese subring $\bigoplus_{i = 0}^{\infty} A_{ni}$ for some $n\geq 2$.
If
$B$ is prime and $A_B$ (or $_B A$) is finitely generated, then $B$ is
projectively simple.
\end{enumerate}
\end{lemma}

\begin{proof}
(a) Suppose that $A$ is not projectively simple. Let $0 \neq I$ be
an ideal of $A$ with $\dim_k A/I = \infty$.  Then $J = I \cap B$ is a
nonzero ideal of $B$ with $\dim_k B/J = \infty$ and $B$ is not
projectively simple.

Conversely, suppose that $A$ is projectively simple and let $0 \neq J$
be an ideal of $B$.  Since $A$ is projectively simple,
$\dim_k A/(A_{\geq n}JA_{\geq n}) < \infty$. Also $A_{\geq n}J A_{\geq n}
\subseteq BJB = J$ for some $n$, so $\dim_k B/J < \infty$ and $B$ is
projectively simple.

(b) Let $I$ be a nonzero graded ideal of $A$. We want to show that
$I\cap B$ is nonzero. Since $A$ is Goldie prime, $I$ contains a
homogeneous regular element of $A$ and thus
$$ \Kdim (A/I)_B < \Kdim A_B =\Kdim B.$$
Hence the map $B\to A/I$ cannot be injective. Thus $B\cap I\neq 0$ and
$B/(B\cap I)$ is finite dimensional because $B$ is projectively simple.
Now $A/I$ is finitely generated over $B/(B\cap I)$, so it is also finite
dimensional.

(c) We think of $B$ as a subring of $A$ which is zero except
in degrees which are multiples of $n$.  By \cite[5.10(1)]{AZ1}, $B$ is
noetherian.
Let $J$ be any nonzero right ideal of $A$.  Then we claim that
$J \cap B \neq 0$.  Suppose this is not true; then if $0 \neq x \in J$ is
any homogeneous element, then $x^n \in J \cap B = 0$.  Thus
$J$ is a right nil ideal, and so $J$ is a nilpotent ideal since $A$ is
noetherian \cite[2.3.7]{MR}.  Then
since $A$ is prime, $J = 0$, a contradiction.

Now let $I$ be a nonzero graded ideal of $B$.  Since $B$ is prime, $I$
contains
a homogeneous regular element $x$.  If $J = \rann_A x$, then $\rann_B x =
J \cap B = 0$
since $x$ is regular in $B$; then by the claim above, $J = 0$ and $x$ is
regular in $A$.
Then $A/IA$ is a noetherian $(B, A)$-bimodule which is Goldie torsion as a
right $A$-module.
By Lemma \ref{lem1.7}(a), $\dim_k A/IA <\infty$.  Since $I = IA \cap B$,
we conclude that $\dim_k B/I<\infty$.
\end{proof}

The result of Lemma \ref{lem1.8}(c) is false without the prime
hypothesis on $B$, as is clear from the following example.

\begin{example}
\label{ex1.9}
Let $A=k\langle x,y\rangle /(x^2, y^2)$. Then $A$ is a PI prime ring of
GK-dimension one, so $A$ is projectively simple by Lemma \ref{lem1.2}(c).
The Veronese subring $A^{(2)}$ is isomorphic to $k[u,v]/(uv)$ where $u=xy$
and $v=yx$. Hence
$A^{(2)}$ is semiprime, but not prime, so it cannot be projectively
simple.
\end{example}

\section{Twisted homogeneous coordinate rings}
\label{xxsec2}

Starting now we work towards the goal of producing some interesting
explicit examples of projectively simple rings. For this we
take Example~\ref{ex1.3} as our model; we expect to find other
projectively simple rings by looking at the class of twisted
homogeneous coordinate rings, which we define and study in this
section.  We will also prove Theorem~\ref{thm0.4}, which will show
that, under certain hypotheses, twisted homogeneous coordinate rings
really are the only examples of projectively simple rings.

Assume throughout this section that $k$ is an algebraically closed
field. Let $X$ be a commutative projective scheme, $\sigma$ an
automorphism of $X$ and $\cal{L}$ an invertible sheaf on $X$. For
any sheaf $\cal{F}$ on $X$, we use the notation $\cal{F}^{\sigma}$
for the pullback $\sigma^{*}(\cal{F})$. Now set
$$\cal{L}_n = \cal{L} \otimes \cal{L}^{\sigma} \otimes \dots
\otimes \cal{L}^{\sigma^{n-1}}$$
for all $n \geq 0$. The \emph{twisted homogeneous coordinate ring} (or
{\it twisted ring} for short) $B = B(X, \cal{L}, \sigma)$ is defined to be
the graded vector space $\bigoplus_{n = 0}^{\infty} \HB^0(\cal{L}_n)$,
with the multiplication rule
$$f g =  f \otimes g^{\sigma^m}\quad \text{ for }\quad f\in B_m,\ g \in
B_n.$$
For more details about this construction see \cite{AV} and \cite{Ke}.

The sheaf $\cal{L}$ is called \emph{$\sigma$-ample} if for any coherent
sheaf
$\mc{F}$ on $X$, $\HB^i(\mc{F} \otimes \mc{L}_n) = 0$ for all $i > 0$ and
$n \gg 0$.
This reduces to the usual notion of ampleness in the commutative case when
$\sigma$
is the identity.  In case $\mc{L}$ is $\sigma$-ample, the ring $B = B(X,
\mc{L},
\sigma)$ is noetherian and there are many nice relationships between the
properties
of $B$ and the properties of $X$. For example, in this case there is an
equivalence
of categories $\Qgr B \simeq \Qch X$, where $\Qch X$ is the category of
quasi-coherent sheaves on $X$ and $\Proj B = (\Qgr B, \pi B)$ is the
noncommutative projective scheme
associated to $B$ (see \S\ref{xxsec10}).


We are now ready to prove Proposition~\ref{prop0.1} from 
the introduction. This result gives a simple geometric criterion for 
projective simplicity of the ring $B(X, \mc{L}, \sigma)$. (For 
the reader's convenience we restate it as
Proposition~\ref{prop2.2} below.) The answer involves the following 
geometric notion.

\begin{definition}
\label{def2.1} Let $\sigma$ be an automorphism of a projective scheme $X$.
Then we
call $\sigma$ \emph{wild} if $\sigma(Y)\neq Y$ for every nonempty reduced
closed
subscheme $Y\subsetneq X$.
\end{definition}

\begin{proposition}
\label{prop2.2}
Let $B = B(X, \mc{L}, \sigma)$, where $\mc{L}$ is $\sigma$-ample.
Then $B$ is projectively simple if and only if $\sigma$ is a wild
automorphism of $X$.
\end{proposition}

\begin{proof} Call two right ideals $J, J'$ of $B$ equivalent if
$J_{\geq n}= J'_{\geq n}$ for some $n \geq 0$.  By \cite[4.4]{AS1},
the mapping $\mc{I} \mapsto \bigoplus_{n \geq 0} \HB^0(\mc{I}
\otimes \mc{L}_n)$ gives a bijective correspondence between ideal
sheaves $\mc{I}$ of $X$ which are $\sigma$-invariant (in other
words $\mc{I}^{\sigma} \cong \mc{I}$), and equivalence classes of
two-sided ideals of $B$.

Note that if $\mc{I}$ is a $\sigma$-invariant ideal sheaf which defines
a closed subscheme $Y$ of $X$, then the ideal sheaf $\mc{I}'$ defining
the reduction $Y_{red}$ of $Y$ is also $\sigma$-invariant.  Moreover, if
$Y$ is already reduced, then its defining ideal sheaf $\mc{I}$ is
$\sigma$-invariant if and only if $\sigma(Y) = Y$. Thus $B$ is
projectively simple if and only if the only reduced subschemes $Y$
of $X$ with $\sigma(Y) = Y$ are $\emptyset$ and $X$.
\end{proof}

Next, we want to prove Theorem~\ref{thm0.4}, which allows us to say
that a projectively simple ring 
that, under certain assumptions, every a projectively simple ring is
a twisted homogeneous coordinate ring in large degree.  Our proof
will largely rely on results from the literature. For the reader's
convenience, we will summarize the main ideas involved,
but without full details.

Let $A$ be a connected finitely $\mb{N}$-graded $k$-algebra. An
$R$-\emph{point module} for $A$, where $R$ is a commutative
$k$-algebra, is a cyclic graded right $A \otimes_k R$-module $M$,
generated in
degree $0$, such that
$M_0 \cong R$ and $M_n$ is a locally free $R$-module of rank $1$ for
each $n \geq 0$. The algebra $A$ is called \emph{strongly noetherian}
if for every commutative noetherian $k$-algebra $C$, the ring
$A \otimes_k C$ is noetherian.   If $A$ is strongly noetherian, then
by the work of Artin and Zhang \cite{AZ2} there exists a projective scheme
$X$ which parametrizes the point modules for $A$.  More exactly, if
we associate to each commutative $k$-algebra $R$ the  set of
all isomorphism classes of $R$-point modules for $A$, then this rule
defines a
functor $\emph{Rings} \to \emph{Sets}$ and $X$ \emph{represents} this
functor.

Now let $A$ be strongly noetherian and generated in degree $1$.  For each
base ring
$R$ one gets a map from the set of $R$-point modules to itself which is
defined by
the rule $M \mapsto M_{\geq 1}(1)$, and these maps induce an automorphism
$\sigma$ of
the representing scheme $X$ \cite[Proposition 10.2]{KRS}.  Then using the
same
methods as in \cite{ATV}, there is a natural way to construct a ring
homomorphism
$\phi: A \to B(X, \mc{L}, \sigma)$ where $\mc{L}$ is a very ample
invertible sheaf on
$X$. In \cite{RZ}, the authors prove the following additional facts about
such maps
$\phi$:

\begin{proposition} \cite[Theorem 1.1]{RZ}
\label{prop2.4}
Let $A$ be a strongly noetherian connected graded ring generated in degree
1.  Let $X$ be the scheme parametrizing the point modules, and let
$\sigma$
be the automorphism of $X$ induced by the map $M \mapsto M_{\geq 1}(1)$ on
point modules.
Then there is a graded ring homomorphism $\phi: A\to B(X, {\cal L},
\sigma)$ which is surjective in
large degree, and such that $\mc{L}$ is $\sigma$-ample.
\end{proposition}

Given the preceding proposition, it is now easy to 
prove Theorem~\ref{thm0.4} from the introduction, 
which we restate below.

\begin{theorem}
\label{thm2.5}
Let $k$ be an algebraically closed field and let $A$ be a
projectively simple noetherian $k$-algebra. Suppose that $A$
is strongly noetherian, generated in degree 1 and has a point module.
Then there is an injective homomorphism $A \hookrightarrow B$
of graded algebras, such that $\dim_k \,  B/A < \infty$
and $B = B(X,{\mathcal L},\sigma)$ is a projectively simple
twisted homogeneous coordinate ring for some smooth projective
variety $X$ with a wild automorphism $\sigma$ and
a $\sigma$-ample line bundle ${\mathcal L}$. 
\end{theorem}

\begin{proof}
Let $\phi$ be the ring homomorphism of Proposition~\ref{prop2.4} above.
By that proposition we have $\im \phi \supset B_{\geq m}$ for all
$m \gg 0$.  Since we
assumed that $A$ has a point module, the scheme $X$ representing the point
modules for $A$
must be nonempty.  Then since $\mc{L}$ is $\sigma$-ample and $X$ is
non-empty, the
sheaf $\mc{L}_n$ has some nonzero global section for all $n \gg 0$
\cite[Proposition 2.3]{Ke}, and so $B_n \neq 0$ for all $n \gg 0$.
Then $\im \phi$ is not finite dimensional, so $\ker \phi$ is an ideal of
$A$ with $\dim_k A/\ker \phi = \infty$.  Since $A$ is projectively simple,
$\ker \phi = 0$ and thus also $\dim_k B/A < \infty$.  The ring $A$ is also
prime by Lemma~\ref{lem1.2}(b), and $B$ must be prime as well.  Then $B$
is projectively simple by Lemma~\ref{lem1.8}(a).  The automorphism
$\sigma$ is wild
by Proposition~\ref{prop2.2}.
\end{proof}

Because of the preceding result, we will focus the majority of our
attention on twisted homogeneous coordinate rings in the remaining
sections of the
paper.  We are very curious, though, what kinds of more general
projectively simple rings may
appear if the hypotheses of Theorem~\ref{thm2.5} are relaxed.
We wonder how stringent the hypothesis that $A$ has a point module
is; do there exist, for example, any projectively simple rings
which have no finitely generated modules of GK-dimension one?

We can offer an example of a projectively simple ring which is
not strongly noetherian.  The example is not too far away from the
strongly noetherian case, since
it has a closely related overring which is strongly noetherian and
projectively simple.  It
also has a large supply of point modules.  Again, it would be interesting
to know if there
are non-strongly noetherian
examples of projectively simple rings which are significantly different
from this one.

\begin{example}
\label{ex2.6}
Let $k$ be algebraically closed of characteristic zero. Let $B =
B(X, \mc{L}, \sigma)$ where $X$ is integral (so $B$ is a domain), and
$\mc{L}$ is $\sigma$-ample and very ample. We construct
special subrings of $B$ which are ``Naive noncommutative blowups''
as studied in \cite{KRS}. Let $c \in X$ be a smooth closed point,
with corresponding ideal sheaf $\mc{I}$.  Setting $\mc{L}_n = \mc{L}
\otimes \mc{L}^{\sigma} \otimes \dots \otimes \mc{L}^{\sigma^{n-1}}$
and $\mc{I}_n = \mc{I} \mc{I}^{\sigma} \dots
\mc{I}^{\sigma^{n-1}}$ for each $n \geq 0$, we define the ring
$$ R = R(X, \mc{L}, \sigma, c)  = \bigoplus_{n \geq 0} \HB^0
(\mc{I}_n \otimes \mc{L}_n) \subset B(X, \mc{L}, \sigma).$$

Now assume in addition that $X$ is an abelian variety, and that the
point $c \in X$ \emph{generates} $X$ in the sense that the subgroup
$\mc{C} = \{ nc | n \in \mb{Z} \}$ of $X$ is a Zariski dense set.
Let $\sigma$ be the translation automorphism $x \mapsto x + c$
of $X$.  In Theorem~\ref{thm6.4} below we will prove that $\sigma$
is wild, and so by Proposition~\ref{prop2.2}, $B = B(X, \mc{L}, \sigma)$
is projectively simple.  By \cite[Theorem 7]{CS}, it follows
that the set $\mc{C}$ is \emph{critically dense}, which means that every
infinite subset of it is dense.  Then $R$ is a noetherian domain by
\cite[Theorem 4.1]{KRS}.

Now let $0 \neq P$ be a homogeneous prime ideal of $R$.  It follows
from \cite[Theorem 4.1 and Lemma 6.1]{KRS} that the $R$-module
$(R/P)$ has a compatible $B$-structure in large degree. More exactly,
there is some $n \geq 0$ and a left $B$-module $M$ such that
$_R M \cong {}_R(R/P)_{\geq n}$. Now either $P = R_{\geq 1}$, or else
$P = \ann_R (R/P)_{\geq n} = (\ann_B M) \cap R$.  In either case,
$P$ is the intersection with $R$ of an ideal $Q$ of $B$.  Since $B$
is projectively simple, $\dim_k B/Q < \infty$, and so $\dim_k R/P
< \infty$.  Thus $R$ is projectively simple.  However, $R$ is not
strongly noetherian when $\dim X\geq 2$ \cite[Theorem 9.2]{KRS}.
\end{example}

\section{Wild automorphisms and algebraic group actions}
\label{xxsec3}

As we saw in the preceding section, to study projectively simple
twisted homogeneous coordinate rings we need to better
understand wild automorphisms.  This will be the subject
of sections~\ref{xxsec3} -- \ref{sec6}. We begin with
several simple observations.

\begin{lemma}
\label{lem2.3}
Let $X$ be a projective scheme. Let $\sigma$ be an automorphism of $X$. 

\begin{enumerate}
\item
If $\sigma$ is wild, then $X$ is reduced.

\item
If $\sigma$ is wild, then $X$ is smooth.

\item
Assume that $X$ is reduced with irreducible components
$X_1, X_2, \dots , X_m$.  Then $\sigma$ is wild 
if and only if the permutation of the $X_i$ induced by $\sigma$ is
a single $m$-cycle, and $\sigma^m$ restricts to a wild automorphism
of each $X_i$. Moreover, if $\sigma$ is wild then $X$ is a disjoint 
union of $X_1,\cdots, X_m$.

\item If $X$ is integral, then $\sigma$ is wild if
and only if $\sigma^n$ is wild for every $n \geq 1$.
\end{enumerate}
\end{lemma}

\begin{proof} (a) $X_{\text{red}}$ is a non-trivial subscheme
of $X$ preserved by $\sigma$; hence, $X = X_{\text{red}}$.

\smallskip
(b) By part (a), $X$ is reduced. Let $Y$ be the singular locus of $X$.
Then $Y$ is
closed, $\sigma$-invariant and $Y \neq X$. Since $\sigma$ is wild, $Y$ is
empty.

\smallskip
(c) The orbit of each component is preserved by $\sigma$; hence, there
can only be one such orbit. Since the subscheme
\[ Y = \bigcup_{i \ne j} \, X_i \cap X_j \subsetneq X \]
is $\sigma$-invariant, it has to be empty, i.e., $X_1, \dots, X_m$ are
disjoint.  If
$\sigma^m(Z) = Z$ for some subscheme $Z \subsetneq X_i$, then $\sigma(Z')
= Z'$ where
$Z' = \bigcup_{j = 0}^{m-1} \sigma^j(Z) \subsetneq X$; thus $Z =
\emptyset$ and so
$\sigma^m \vert_{X_i}$ is wild.

\smallskip
(d) If $\sigma$ preserves a subscheme $Y \subset X$, then so does
$\sigma^n$.
Conversely, if $\sigma^n$ preserves $Y \subsetneq X$ then $\sigma$
preserves
$\bigcup_{j = 0}^{n-1} \sigma^j(Y) \subsetneq X$.
\end{proof}

In view of Lemma~\ref{lem2.3}, it is clear that we lose nothing essential in
our understanding of wild automorphisms by restricting to the case where $X$
is integral and smooth, and we will usually do so.

The rest of this section will be devoted to proving
Conjecture~\ref{con0.3} under the additional assumption that $\sigma$
can be extended to a finite-dimensional family of automorphisms
(or even of birational isomorphisms) of $X$.  Before proceeding,
let us make a few remarks about our terminology. In order to avoid
arithmetic complications (in what is already a difficult geometric
problem),
we shall assume throughout that $k$ is an algebraically closed field
of characteristic zero.  For convenience, we often use the language of
varieties rather than
schemes; for us, a \emph{variety} is a (not necessarily irreducible)
reduced
separated scheme of finite type over $k$.  For a variety $X$
the notation $p \in X$ always means that $p$ is a closed point of $X$.
By an \emph{algebraic group} we mean a variety with a compatible group
structure.  Note  that we do not assume that algebraic groups are linear,
but by definition they must have finitely many components. If $G$ is an
algebraic group, then for each $h \in X$ one has a
\emph{translation automorphism} $T_h$
defined by $g \mapsto hg$.  An \emph{abelian variety} is a complete
irreducible variety $X$ over $k$ which is also an algebraic group.

\begin{proposition}
\label{thm5.3}
 Let $X$ be an irreducible projective variety. Suppose that an algebraic
 group $G \subset \Aut(X)$ acts regularly on $X$ such that $\sigma\in G$
 acts by a wild
 automorphism. Then $X$ is an abelian variety and some power of $\sigma$
 is a translation automorphism when a group structure in $X$ is chosen
 properly.
\end{proposition}

\begin{proof}
First we make some reductions. Let $H$ be the closure in $G$ of the
subgroup $\{ \sigma^i \}_{i \in \mb{Z}}$. Then $H$ is an algebraic group.
Since $H$ is a closure of an abelian group, $H$ is abelian.  Without loss
of generality, we may assume from now on that $G=H$, so in particular
$G$ is abelian.  Since $G$ has finitely many components, $\sigma^n$ must
lie in the connected component $G_e$ containing the identity element of
$G$, for some $n \geq 1$. Since $X$ is irreducible, $\sigma^n$ is also
wild (Lemma~\ref{lem2.3}(d)).  Thus we may replace
$\sigma$ by $\sigma^n$ and $G$ by $G_e$, and
so we may assume that $G$ is irreducible.

Choose any $x \in X$, and let $Y$ be the closure of $Gx$ in $X$. Since
$\sigma \in G$, $\sigma(Gx) \subset Gx$ and hence $\sigma(Y) \subset Y$.
Since $\sigma$ is wild, $Y=X$. Now the rule $g \mapsto gx$ defines a
morphism $f:G\to X$. By Chevalley's theorem \cite[Exercise II.3.19]{Hart},
the image $Gx$ of the map $f$ must be constructible, and so it must
contain
an open subset $U$ of its closure $Y=X$. Hence $Gx=GU=\bigcup_{g\in G} gU$
is dense and open in $X$.  Since $Gx$ is $\sigma$-stable, $Z:=X-Gx$ is a
$\sigma$-stable closed subvariety. Since $\sigma$ is wild, we conclude
that $Z$ is empty and thus $X=Gx$.

Next, let $G_0$ be the stablizer of $x$. Since $G$ is abelian, $G_0$ is
the stablizer of every point in $Gx = X$. Since $G$ is a subgroup of
$\Aut(X)$ and $X$ is a variety over an algebraically closed field,
automorphisms
of $X$ are determined by their actions on closed points and so $G_0$
is trivial.  Hence the morphism $f: G\to X$ is bijective.

Since our standing assumption is that the base field $k$ has
characteristic $0$, it now follows
from \cite[Theorem 4.6]{Hum} that $f$ is birational.  Let $V$ be the
largest
open set of $X$ such that $f \vert_{f^{-1}(V)}: f^{-1}(V) \to V$ is an
isomorphism.  Then $V$ is $\sigma$-invariant, and since $\sigma$ is wild,
we see that $V = X$ and $f$ is an isomorphism.  Therefore the group
structure of $G$ may be transferred to $X$ via $f$, $G$ must be
projective,
and both $X$ and $G$ are abelian varieties.  The isomorphism $f$
transforms
the translation automorphism $T_\sigma$ of $G$ to the translation
automorphism $T_{f(\sigma)} = \sigma$ of $X$.  Since we replaced $\sigma$
by $\sigma^n$ during the proof, we see that some power of the original
$\sigma$ is a translation automorphism of $X$ as required.
\end{proof}

The previous theorem can be generalized to the case of birational actions.
Let $\Bir(X)$ be the group of all birational maps over $k$ from $X$ to
itself.

\begin{corollary}
\label{cor5.4}
Let $X$ be an irreducible projective variety admitting a wild automorphism
$\sigma$.  Suppose that $P\subset \Bir(X)$ is an algebraic group that
acts birationally on $X$ such that $\sigma\in P$. Then $X$ is an abelian
variety and some power of $\sigma$ is a translation after a group
structure in $X$ is chosen properly.
\end{corollary}

\begin{proof} Let $G=\{\tau\;|\; \sigma \tau=\tau\sigma\}$ be the
centralizer of $\sigma$ in $P$. Since $P$ is an algebraic group, so is
$G$.
Clearly $\sigma\in G$. We want to show that $G\subset \Aut(X)$. For any
$\tau\in G$, let $Y_{\tau}$ be the indeterminacy locus of the rational map
$\tau \colon X \dasharrow X$.
Since $\sigma\tau=\tau\sigma$, $Y_\tau$ is $\sigma$-stable. Since
$\sigma$ is wild,
$Y_\tau$ is empty, whence $\tau$ is regular. Similarly, $\tau^{-1}$ is
regular. Therefore $G\subset \Aut(X)$. The assertion now follows from
Theorem \ref{thm5.3}.
\end{proof}

In Remark \ref{rem7.5} below, we will see that the conclusions of
Theorem \ref{thm5.3} and Corollary \ref{cor5.4} may be strengthened to
say that $\sigma$ is, in fact, itself a translation automorphism.

\section{Wild automorphisms and numerical invariants}
\label{sec5b}

In this section we will show that the existence of a wild automorphism
imposes strong restrictions on two important numerical invariants
of a projective variety, the Kodaira dimension and the Euler
characteristic. We will continue to assume that the base field $k$
is algebraically closed and of characteristic zero.

We begin by considering the Kodaira dimension.  Suppose that $X$
is a smooth projective variety, and let $\omega_X$ be the
canonical sheaf.  For a sheaf $\mc{F}$ on $X$ we
write $\Gamma(X, \mc{F})$ for the global sections of $\mc{F}$. The
\emph{Kodaira dimension} (or \emph{canonical dimension}) of $X$, denoted
by $\kappa(X)$, is defined to be the transcendence degree of the canonical
ring $\bigoplus_{i\geq 0} \Gamma(X, \omega_X^{\otimes i})$ minus $1$.
Similarly, the {\it anti-canonical dimension} $\bar{\kappa}(X)$ of $X$
is the transcendence degree of the anti-canonical ring
$\bigoplus_{i\geq 0} \Gamma(X, \omega_X^{\otimes -i})$ minus $1$.
We will
see next that the existence of a wild automorphism limits
the possible values of $\kappa$ and $\bar{\kappa}$.

Let ${\cal L}$ be an invertible sheaf on a projective variety $X$ and
let $\sigma $ be an automorphism of $X$. A
$\sigma$-linearization of ${\cal L}$ is an isomorphism
$$ p \colon \sigma^{\ast}({\cal L}) \stackrel{\simeq}{\longrightarrow}
{\cal L} \, . $$
Informally speaking, we can think of $p$ as a way to lift $\sigma$ from
$X$ to ${\cal L}$. If $p$ is chosen, we shall say that ${\cal L}$ is
$\sigma$-linearized and identify $\sigma^{\ast}({\cal L})$ with ${\cal L}$
via $p$.

\begin{proposition}
\label{prop5.5}
Suppose an irreducible projective variety $X$ admits a wild automorphism
$\sigma$.  If ${\cal L}$ is a $\sigma$-linearized invertible sheaf on $X$
then either
\begin{enumerate}
\item
${\cal L}$ is isomorphic to ${\cal O}_X$, or
\item
$\Gamma(X, \mc{L}) = 0 = \Gamma(X, \mc{L}^{-1})$.
\end{enumerate}
\end{proposition}

\begin{proof} Suppose that (b) fails, in other words that either $\mc{L}$
or $\mc{L}^{-1}$ has a nonzero global section.
Note that a $\sigma$-linearization of ${\cal L}$
induces a $\sigma$-linearization of $\mc{L}^{-1}$. Thus
after switching the roles of $\mc{L}$ and $\mc{L}^{-1}$
if necessary, we may assume that $V := \Gamma(X, \mc{L}) \neq (0)$.

Since ${\cal L}$ is linearized, pullback by $\sigma$
induces an automorphism $\wt{\sigma} = p \circ \sigma^{\ast}$ of $V$.
Since $k$ is algebraically closed, $V$ must have
an eigenvector $f$ for the action of $\wt{\sigma}$.
The vanishing set $Z(f) \subsetneq X$ of the global section $f$
is then fixed by $\sigma$.  Since $\sigma$ is wild,
$Z(f) = \emptyset$. This implies that ${\cal L}$ is
generated by the global section $f$, i.e.,
the map $i : \mc{O}_X \to \mc{L}$ sending $1$ to $f$ is an isomorphism.
\end{proof}

\begin{corollary}
\label{cor5.6}
Suppose an irreducible projective variety $X$ of dimension $d \ge 1$
admits a wild automorphism.
Then either (i) $\kappa(X) = \bar{\kappa}(X) =  -1$ or (ii)
$\omega_X^{\otimes n} \cong \mc{O}_X$ for some $n \ge 1$. (In the latter
case $\kappa(X) = \bar{\kappa}(X) = 0$.)
In particular, $X$ cannot be a Fano variety or a variety
of general type.
\end{corollary}

\begin{proof} Apply Proposition \ref{prop5.5} with ${\cal L} =
\omega_X^{\otimes n}$. If (a) holds for some $n \ge 1$
then $\omega_X^{\otimes n} \cong \mc{O}_X$; if instead (b)
holds for all $n \ge 1$ then $\kappa(X) = \bar{\kappa}(X) =  -1$.
The last assertion is an immediate consequence of the first:
if $X$ is Fano then $\bar{\kappa}(X) = d \ge 1$ and if $X$
is of general type then $\kappa(X) = d \ge 1$.
\end{proof}

\begin{corollary}
\label{wild auts of subvarieties} Let $X$ be an abelian variety, and
suppose that $Y
\subseteq X$ is an irreducible subvariety of $X$ such that $Y$ admits a
wild
automorphism.  Then $Y$ is a translate of an abelian subvariety of $X$.
\end{corollary}
\begin{proof}
By Corollary~\ref{cor5.6}, since $Y$ has a wild automorphism it follows
that
$\kappa(Y) \leq 0$.  Then a theorem of Ueno, which is proved for the case
$k =
\mb{C}$ in \cite[10.1 and 10.3]{Ue}, and for a general algebraically closed
field in
\cite[3.7]{Mo}, states that $Y$ must be a translate of an abelian
subvariety of $X$.
\end{proof}

In the rest of this section we work out what the existence
of a wild automorphism tells us about the Euler
characteristic (or the arithmetic genus) of a variety.
As a consequence, we will see that rationally connected varieties
have no wild automorphisms.

Given a regular map $f \colon X \lra X$, we define the {\em algebraic
Lefschetz number} $L(f, X)$ by the formula
$$ L(f,X) =\sum_{q \geq 0} (-1)^q \text{Trace} \,
\bigl(f^* \colon \HB^q(X,\cal{O}_X)
\lra \HB^q(X, \cal{O}_X) \bigr) \, . $$
Note that the sum on the right is well-defined because $\HB^q(X,
\cal{O}_X)$
is a finite-dimensional $k$-vector space for every $q \ge 0$
\cite[III.5.2]{Hart} and  $\HB^q(X, \cal{O}_X) = (0)$ for $q > \dim(X)$
\cite[III.2.7]{Hart}. We also define the {\em algebraic Euler
characteristic} $\chi(\cal{O}_X)$ of $X$ as
$$\chi(\cal{O}_X) =
L(\id ,X) =\sum_{q \geq 0} (-1)^q \dim \, \HB^q(X,\cal{O}_X) \, .$$

\begin{proposition}
\label{prop5.8}
Let $X$ be a smooth irreducible projective variety and let $\sigma$ be
an automorphism of $X$.

\smallskip
(a) If $\chi({\cal O}_X) \neq 0$, then some
power of $\sigma$ has a fixed point.

\smallskip
(b) If $\dim(X) \ge 1$ and $\sigma$ is wild then the algebraic
Euler characteristic $\chi(\cal{O}_X) = 0$. Equivalently,
the arithmetic genus $p_a(X) = (-1)^{\dim(X) + 1}$.
\end{proposition}

\begin{proof}
(a) By the ``holomorphic Lefschetz fixed point theorem",
if $\sigma$ has no fixed points then $L(\sigma, X) = 0$.
A proof of this result over $k = \mathbb{C}$
can be found in~\cite[Section III.4, p. 426]{GH} and over an
arbitrary algebraically closed field $k$ of characteristic zero in
\cite{TT}. (Note that the sum in the right hand side of the formula
is empty if there are no fixed points.)

Now assume that $\sigma^n$ has no fixed point for any $n\geq 1$.
Set $W=\bigoplus_{q \geq 0} \HB^q(X,\cal{O}_X)$ and $d = \dim_k(W)$.
Then the linear map $\sigma^* \colon W \lra W$ induced by $\sigma$
satisfies a characteristic equation of the form
$$
(\sigma^*)^d+c_1 (\sigma^*)^{d-1}+\cdots +c_d \, \id_W = 0
$$
for some $c_i \in k$. Since $\sigma^*$ is an automorphism of $W$,
$c_d = (-1)^d \det(\sigma^*) \neq 0$.
By the linearity of the trace, we have
$$
L(\sigma^d,X)+c_1 L(\sigma^{d-1},X)+ \dots + c_{d-1} L(\sigma, X) +
c_d L(\id_X, X)=0.
$$
Then since $L(\sigma^n,X)=0$ for all $n\geq 1$,
$L(\id_X,X)=\chi(\cal{O}_X) =0$.

\smallskip
(b) By Lemma~\ref{lem2.3}(d) no power of $\sigma$ can have
a fixed point.  Hence, part (a) tells us that $\chi(\cal{O}_X) = 0$.
The second assertion is simply a restatement
of the first, since by definition,
$p_a(X) = (-1)^{\dim(X)} (\chi(\cal{O}_X) -1)$ \cite[p. 230]{Hart}.
\end{proof}

\begin{remark}
\label{rem5.10}
If $X$ is an irreducible variety
defined over $k = \mathbb{C}$ which has a wild automorphism, then the same
arguments as in Proposition~\ref{prop5.8} show that $e(X) = 0$,
where $e(X)$ is the usual (topological) Euler characteristic.
\end{remark}

\begin{proposition} \label{prop5.11}
Suppose an irreducible projective variety $X$ of dimension $d \ge 1$ has a
wild
automorphism. Then

\smallskip
(a) $X$ carries a non-vanishing regular differential
$m$-form $\omega$ for some odd $m \ge 1$.  Moreover,
$\omega$ can be chosen so that
\begin{equation} \label{e.eigenvalue}
\text{$\sigma^{\ast} \omega = c \, \omega$ for some $c \in k^*$.}
\end{equation}

(b) $X$ is not rationally connected.

\smallskip
(c) $X$ is not unirational.
\end{proposition}

Recall that an algebraic variety
$X$ is called {\em rationally connected} if there is a family of
irreducible rational curves in $X$ such that two points
in general position can be connected by a curve from
this family; see~\cite[Section IV.3]{kollar}.
$X$ is called \emph{unirational} if it admits a dominant
rational map $\mb{P}^n \dasharrow X$ for some $n \ge 1$.

\begin{proof}
(a) Let $h^{n, m} = h^{n, m}(X) = \dim \, \HB^n(X, \bigwedge^m \Omega_X)$
be the Hodge numbers of $X$. Here $\Omega_X$ is the sheaf of differential
$1$-forms on $X$ and $\bigwedge^m \Omega_X$ is the sheaf
of differential $m$-forms.

We claim that $h^{0, m} > 0$ for some odd $m \ge 1$.
Since $h^{n, m} = h^{m, n}$ (see, for example,
\cite[pp. 54-55]{danilov}),
we have \[ p_a(X) = \sum_{i = 0}^{d-1} (-1)^i h^{d-i, 0} =
\sum_{m = 1}^{d} (-1)^{d-m} h^{0, m} \, . \]
By Proposition \ref{prop5.8}, $p_a(X) = (-1)^{d+1}$, i.e.,
\begin{equation}
\label{e.hodge}
\sum_{m = 1}^{d} (-1)^{m + 1} h^{0, m} = 1 \, .
\end{equation}
If $h^{0, m} = 0$ for every odd $m \ge
1$, then every term in the left hand side of~\eqref{e.hodge} is
non-positive, a
contradiction.

Thus for some odd $m$ we have $V := \HB^0(X, \bigwedge^m \Omega_X) \ne
(0)$.  Now
we can choose $\omega$ to be an eigenvector for the linear
automorphism $\sigma^*$ of the finite-dimensional $k$-vector space $V$.
In other words,~\eqref{e.eigenvalue} holds.  Finally, the
vanishing locus $Y$ of $\omega$ is a closed
$\sigma$-invariant subvariety of $X$. Since $\sigma$ is
wild, $Y = \emptyset$, i.e., $\omega$ is a non-vanishing differential
form, as claimed.

\smallskip
(b) follows from (a). Indeed, if $X$ is rationally connected
then $h^{0, m}(X) = 0$ for
all $m \ge 0$; see~\cite[Corollary IV.3.8]{kollar}.

\smallskip
(c) follows from (b) because a unirational variety is rationally
connected;
see~\cite[Example IV.3.2.6.2]{kollar} or~\cite[Lemma 3.4.1]{lbr}.
\end{proof}

\section{Wild automorphisms and the Albanese map}
\label{sec5c}

Let $X$ be an irreducible variety.  Then associated to $X$ is the
\emph{Albanese
variety} $\Alb(X)$, which is an abelian variety, and the \emph{Albanese
map} $f: X
\to \Alb(X)$.  These constructions have the following properties:
$\Alb(X)$ is
generated as an algebraic group by the image of $f$ (in the sense of
Definition~\ref{def.generate} below), and given any other regular map $g:
X \to Y$
where $Y$ is an abelian variety, there is a regular map $h: \Alb(X) \to Y$
such that
$hf = g$.  The dimension of $\Alb(X)$ is called the~\emph{irregularity} of
$X$ and is
denoted by $q(X)$ (or just $q$ if the reference to $X$ is clear from the
context).
Equivalently, $q(X) = \dim_k \HB^1(X, \mc{O}_X)$; see for example
\cite[Section
I.6]{BH} or \cite[Lemma 9.22]{Ue}.

We shall now see that if an irreducible variety has a wild automorphism
then its Albanese map must have various special properties.
\begin{proposition}
\label{prop5.7} Let $X$ be an irreducible projective variety of dimension
$d$ and irregularity $q$
and let $\pi \colon X \lra \Alb(X)$ be the Albanese map. Suppose
that $\sigma$ is a wild automorphism of $X$.
By the universal property of the Albanese map,
$\sigma$ induces an automorphism of $\Alb(X)$, which
we will denote by $\overline{\sigma}$.

\begin{enumerate}
\item
$\pi$ is surjective.
\item
$\overline{\sigma}$ is a wild automorphism of $\Alb(X)$.
\item
$\pi$ is smooth.
\item
The fiber $X_t = \pi^{-1}(t)$ is a smooth irreducible variety
of dimension $d - q$ for every $t \in \Alb(X)$.
\item
$q \le d$ and if $q = d$ then $X$ is an abelian variety.
\end{enumerate}
\end{proposition}

Note that by~\cite[Theorem 8.1]{Mo}, part (a) holds for
every variety $X$ with $\kappa(X) = 0$.  On the other hand,
if $\kappa(X) = -1$ and we do not assume that $X$ has a wild
automorphism, then part (a) may fail; for example, the Albanese
map for the ruled surface $X \lra C$ is surjective if and
only if $C$ is a curve of genus $\le 1$; see~\cite[p. 554]{GH}.

\begin{proof}
(a) Let $\overline{X} = \pi(X)$.  We claim that
\begin{equation} \label{e.alb-a}
\text{ $\overline{\sigma}$ restricts to a wild automorphism
of $\overline{X}$.}
\end{equation}
Indeed, assume the contrary, say
$\overline{\sigma}(\overline{Y}) \subset \overline{Y}$ for $\emptyset
\subsetneq
\overline{Y} \subsetneq \overline{X}$. Then setting $Y =
\pi^{-1}(\overline{Y})$, we see that
$\emptyset \subsetneq \sigma(Y) \subseteq Y$, a contradiction.
Now by Corollary \ref{wild auts of subvarieties},
$\overline{X}$ is a translate
of an abelian subvariety in $\Alb(X)$.  By the definition
of the Albanese map, this implies that
$\overline{X} = \Alb(X)$; see \cite[II.3]{La1}.

\smallskip
(b) This follows from part (a) and~\eqref{e.alb-a}.

\smallskip
(c) Since $\ch k = 0$ and we know that $X$ must be nonsingular
[Lemma~\ref{lem2.3}(b)],
by generic smoothness there exists a non-empty Zariski open 
subset $U \subseteq \Alb(X)$ such
that $\pi$ is smooth over $U$ \cite[Corollary III.10.7]{Hart}.  Then
clearly $\pi$ is
smooth over the $\overline{\sigma}$-invariant open subset
$$ W = \bigcup_{i \in \mathbb{Z}}\; \overline{\sigma}^i(U) $$
of $\Alb(X)$. Since $\overline{\sigma}$ is wild, $W = \Alb(X)$.

\smallskip
(d) The fact that each $X_t$ is smooth of dimension $d-q$
is immediate from part (c). To show that $X_t$ is irreducible,
consider the Stein decomposition
\[ \pi \colon X \stackrel{\alpha}{\lra} X' \stackrel{\beta}{\lra} \Alb(X)
\, , \]
where $\alpha$ has connected fibers and $\beta$ is finite.  The variety
$X'$ may be defined explicitly as $\mathbf{Spec}\ \pi_{*}(\mc{O}_X)$
\cite[III.1.5]{Hart}; in particular,
$X'$ is irreducible, since it is covered by the spectra of domains.
Our goal is to show that $\beta$ is an isomorphism.
If we can prove this, then each $X_t$ is connected and nonsingular, so
must be irreducible.

Since the automorphism $\overline{\sigma}$ of $\Alb(X)$ acts on the
sheaf of
graded algebras $\pi_*(\mc{O}_X)$ over $\Alb(X)$, we get an induced
automorphism
$\wt{\sigma}$ of $X'$.  Since $\beta$ is a finite surjective map and
$\overline{\sigma}$ is
a wild automorphism of $\Alb(X)$, it follows that $\wt{\sigma}$ is a wild
automorphism of $X'$.
Thus $X'$ is nonsingular by Lemma~\ref{lem2.3}(b).
Then as in the proof of part (c), since $\beta$ is smooth over an open set
and $\overline{\sigma}$ is wild, $\beta$
must be a smooth morphism of relative dimension $0$, in other words an
{\'e}tale map.
Now by a theorem of Serre and Lang~\cite[Section IV.18]{Mu},
$X'$ has a structure of an abelian variety such that
$\beta$ is a regular homomorphism. By the universal property
of the Albanese map, $\alpha$ factors through $\pi$.
In other words, $\beta$ has an inverse, so $\beta$
is an isomorphism as desired.

\smallskip
(e) The inequality $q \le d$ is an immediate consequence of part (a).
If $q = d$, then by (c) $\pi$ is a smooth morphism of relative
dimension $0$, so {\'e}tale. Applying \cite[Section IV.18]{Mu} once again, 
we conclude that $X$ has the structure of an abelian variety.
\end{proof}

\section{Wild automorphisms of algebraic surfaces}
\label{sec5d}

In this section we will prove Conjecture~\ref{con0.3} in the case
where $\dim(X) \le 2$. Our proof relies on the classification of
algebraic surfaces. We will start with the most difficult cases,
where $X$ is assumed to be a ruled surface or a hyperelliptic surface,
and defer the rest of the proof until the end of this section.

Recall that an algebraic surface $X$ is called \emph{ruled} if it
admits a surjective morphism
\begin{equation} \label{e.ruled1}
\pi \colon X \to C
\end{equation}
to a smooth curve $C$ whose fiber over each 
point of $C$ is isomorphic to $\mb{P}^1$.  Such
a morphism always has a section
$s \colon C \lra X$; see, e.g., \cite[Corollary IV.6.6.2]{kollar}.

\begin{lemma} \label{lem.ruled}
A (smooth projective) ruled surface cannot have a wild automorphism.
\end{lemma}

\begin{proof} Assume the contrary: $\sigma$ is a wild automorphism of
a ruled surface $X$. Let $\pi: X \lra C$ be as in~\eqref{e.ruled1} and
let $g$ be the genus of $C$.  Since  $p_a(X) = - g$
(see~\cite[Corollary V.2.5]{Hart}),
Proposition~\ref{prop5.8} tells us that $g = 1$, so that $C$
is an elliptic curve.

By~\cite[Corollary V.2.5]{Hart} the irregularity
$q(X) = g = 1$, so the Albanese variety $\Alb(X)$ is
an elliptic curve. We claim that $C$ is, in fact, the
Albanese variety and $\pi$ is the Albanese map for $X$.
Indeed, let $\alpha \colon X \lra \Alb(X)$
be the Albanese map. By the universal property of
the Albanese map, $\pi$ factors through $\alpha$, in other words
\[ \pi \colon X \xrightarrow{\alpha} \Alb(X) \xrightarrow{\beta} C \, , \]
where $\beta$ is a covering map of elliptic curves. Since the fibers
of $\pi$ are connected (each fiber is isomorphic to $\mathbb{P}^1$),
we conclude that $\beta$ is one-to-one, i.e., is a bijective
morphism between smooth curves $\Alb(X)$ and $C$. Consequently,
$\beta$ is an isomorphism; this proves the claim.
(For a different proof of this claim over
$k = \mathbb{C}$, see~\cite[p. 554]{GH}.)

Now let $s \colon C \lra X$ be a section of $\pi$,
with the property that $C_0 = s(C) \subset X$ has the minimal possible
self-intersection number. Following~\cite[Section V.2]{Hart},
we denote this number by $-e$.  The group $\Num(X)$ of divisors
in $X$ up to numerical equivalence is isomorphic to $\mathbb{Z} \oplus
\mathbb{Z}$, and is generated by $C_0$ and $F$, where $F$ is a fiber
of $\pi$; the intersection form on $\Num(X)$ is given
by $C_0 \cdot F = 1$, $F^2 = 0$ and $C_0^2 = - e$ \cite[Proposition
V.2.3]{Hart}.  As we saw in
Proposition~\ref{prop5.7}, $\sigma$ acts on the fibers of $\pi$; since
these fibers are all algebraically
equivalent, we have $\sigma(F) \equiv F$ in $\Num(X)$.
We claim that $C_0$ and $\sigma(C_0)$
are numerically equivalent, that is,
\begin{equation} \label{e.C_0}
\text{$\sigma(C_0) \equiv C_0$ in $\Num(X)$.}
\end{equation}
Indeed, suppose that $\sigma(C_0)\equiv a C_0 + b F$ for some
$a, b \in \mathbb{Z}$. Then $a = \sigma(C_0) \cdot F = 1$, and
since $C_0^2 = \sigma(C_0)^2$, we see that $b = 0$, thus
proving~\eqref{e.C_0}. We now consider three cases.

\smallskip
1. $e > 0$. In view of~\eqref{e.C_0},
$C_0 \cdot \sigma(C_0) = C_0^2 < 0$; consequently, $\sigma(C_0) = C_0$,
contradicting the fact that $\sigma$ is wild.

\smallskip
2. $e = 0$. Since $\sigma$ is wild (and hence, so is $\sigma^2$),
we may assume that $C_0$, $\sigma(C_0)$ and $\sigma^2(C_0)$ are
three distinct curves in $X$.  Since $C_0^2 = 0$,
formula~\eqref{e.C_0} tells us that $C_0$, $\sigma(C_0)$
and $\sigma^2(C_0)$ are mutually disjoint in $X$.
We now appeal to the general fact
that any ruled surface $X \lra C$ with
three mutually disjoint sections $C_1$, $C_2$ and $C_3$
is isomorphic to $\mathbb{P}^1 \times C$ (over $C$). Indeed,
for such $X$ there is a (unique) isomorphism
$\mathbb{P}^1 \times C \lra X$ (over $C$), sending
$\{ 0 \} \times C$, $\{ 1 \} \times C$ and $\{ \infty \} \times C$
to $C_0$, $C_1$ and $C_2$.  Thus we may assume
that $X = \mathbb{P}^1 \times C$, where $C$ is an elliptic curve.
In this case the canonical divisor of $X$ is
$K_X = -2 (\{ \pt \} \times C)$ \cite[Lemma V.2.10]{Hart}, so that
$\bar{\kappa}(X) > 0$, contradicting Corollary~\ref{cor5.6}.

\smallskip
3. $e < 0$. By \cite[Theorems V.2.12 and V.2.15]{Hart} there
is only one surface $X$ in this category, with $e = -1$.
(This surface is often denoted by ${\bf P_1}$.)
By \cite[Theorem 3(4)]{Mar},
the automorphism group $\Aut(X)$ has a normal subgroup
\[ \Delta = \{ \id, \delta_1, \delta_2, \delta_3 \}
\simeq \mathbb{Z}/2 \mathbb{Z} \times \mathbb{Z} / 2 \mathbb{Z} \]
such that $\overline{\delta_1} = \overline{\delta_2} =
\overline{\delta_3} = \id_C$ for $i = 1, 2, 3$.
Here $\overline{\delta}$ is the automorphism of $C$ induced by $\delta$
(as in Proposition~\ref{prop5.7}(b)). Let
\[ Y = X^{\delta_1} \cup X^{\delta_2} \cup X^{\delta_3} \, , \]
where $X^{\delta}$ denotes the fixed point set of $\delta$. Since
$\sigma (X^{\delta}) = X^{\sigma \delta \sigma^{-1}}$, we see
that $\sigma(Y) = Y$. Clearly, $Y \neq X$; it thus remains
to show that $X^{\delta} \neq \emptyset$ for
$\delta = \delta_1, \delta_2, \delta_3$.  Indeed,
for every $p \in C$ we have an automorphism
\[ \delta_{| \, \pi^{-1}(p)} \colon \mathbb{P}^1 \lra \mathbb{P}^1,  \]
where $\mathbb{P}^1$ stands for $\pi^{-1}(p)$.  Since every automorphism
of $\mb{P}^1$ has at least one fixed point, we see that
$X^{\delta} \neq \emptyset$, as claimed.

This completes the proof of Lemma~\ref{lem.ruled}.
\end{proof}

Next we consider hyperelliptic surfaces. Recall that a hyperelliptic
surface $X$ is a surface of the form $(E \times F)/G$, where $E$ and $F$
are elliptic curves, and $G$ is a finite subgroup of $E$
acting on $E \times F$ via
\[ g \cdot (x, y) \mapsto (x + e, \alpha(g) \cdot y) \]
for some $e \in X$ and injective homomorphism $\alpha \colon G
\hookrightarrow
\Aut(F)$. There are only seven possibilities for $E$, $F$, $G$, and
$\alpha$;
see~\cite[Table 1.1]{BM}, \cite[p. 590]{GH} or \cite[p. 37]{bomu}.

\begin{lemma} \label{lem.hyperelliptic}
A hyperelliptic surface cannot have a wild automorphism.
\end{lemma}
\begin{proof}
Suppose that $X$ is a hyperelliptic surface and $\sigma$ is
an automorphism of $X$.
By~\cite[Lemmas 1.2 and 2.1]{BM}, $\sigma$ can be lifted
to an automorphism $\wt{\sigma} \in \Aut(E) \times \Aut(F)$
which normalizes $G$ and such that the diagram
\[ \begin{array}{ccc}
E \times F  & \stackrel{\wt{\sigma}}{\lra} & E \times F \\
  \downarrow &  & \downarrow \\
X = (E \times F)/G & \stackrel{\sigma}{\lra} & X = (E \times F)/G
\end{array} \]
commutes. (In fact, by~\cite[Lemma 3.1]{BM} $\wt{\sigma}$
centralizes $G$, but we shall not use this.)
Consequently, $\sigma$ descends to an automorphism $\sigma_0$
of $F/\alpha(G) \simeq
\mathbb{P}^1$ (cf.~\cite[Theorem 4, pp. 35-36]{bomu}) such that the
extended diagram
\[ \begin{array}{ccc}
E \times F  & \stackrel{\wt{\sigma}}{\lra} & E \times F \\
  \downarrow &  & \downarrow \\
 X = (E \times F)/G & \stackrel{\sigma}{\lra} & X = (E \times F)/G \\
  \downarrow &  & \downarrow \\
 \mathbb{P}^1 \simeq F/ \alpha(G) & \stackrel{\sigma_0}{\lra} &
\mathbb{P}^1 \simeq F/ \alpha(G) \end{array} \] commutes; here the
vertical arrows
are the natural projections. Now recall that every automorphism of
$\mathbb{P}^1$ has
a fixed point. If $x \in F/\alpha(G) \simeq \mathbb{P}^1$ is a fixed point
of
$\sigma_0$ then the preimage of $x$ in $X = (E \times F)/G$ is a
$\sigma$-invariant
curve. This shows that $\sigma$ cannot be wild.

The above argument relied on the results of~\cite{BM} which
are stated under the assumption that $k = \mathbb{C}$.
However, we observe that the the proofs
of Lemmas 1.2 and 2.1 in~\cite{BM} remain valid
over any algebraically closed field $k$ of characteristic $0$.
\end{proof}

\begin{theorem}
\label{thm5.12} Let $X$ be an irreducible projective variety
of dimension $\le 2$ over an algebraically closed
field $k$ of characteristic $0$.  If $X$ admits a wild
automorphism then $X$ is an abelian variety.
\end{theorem}

\begin{proof} By Lemma~\ref{lem2.3}(b), $X$ is smooth.
The case $\dim X = 0$ is trivial. Suppose $\dim(X) = 1$.
Recall that the arithmetic genus $p_a(X)$ and the geometric
genus $p_g(X)$ of a smooth curve $X$ coincide;
see~\cite[Proposition IV.1.1]{Hart}.
Hence, by Proposition~\ref{prop5.8}, $p_g(X) = 1$, i.e.,
$X$ is an elliptic curve or equivalently, a 1-dimensional abelian variety.

 From now on we will assume that $X$ is a smooth surface.
By Proposition~\ref{prop5.8}, $p_a(X) = -1$ and thus the irregularity is
\[ q(X) = p_g(X) - p_a(X) = p_g(X) + 1 \, ; \]
see \cite[Remark 7.12.3]{Hart}. Here
the geometric genus $p_g(X) = \dim \HB^0(X, \omega)$
is either $0$ or $1$; see Proposition~\ref{prop5.5}.
If $p_g(X) = 1$ then $q(X) = 2$, and
$X$ is an abelian variety; see Proposition~\ref{prop5.7}(e).

Now suppose $p_g(X) = 0$ and $q(X) = 1$. To complete the proof,
it is enough to show that this case cannot occur.
Let $\pi \colon X \lra C$ be the Albanese map. Since $q(X) = 1$,
$C$ is an elliptic curve. We claim that $X$ is minimal, i.e.,
$X$ contains no smooth
rational curves $D$ with $D^2 = -1$. Indeed, if such a curve existed
on $X$, we would have $\pi(D) = \{ p \}$, where $p$ is a point of $C$.
In other words, $D$ would be contained in a fiber of $\pi$.
By Proposition~\ref{prop5.7},  $\pi$ is a smooth map with
irreducible fibers.  In particular, $D = \pi^{-1}(p)$. But
then $D^2 = 0$ (because the fibers of $\pi$ are disjoint and
algebraically equivalent). This contradicts our assumption
that $D^2 = -1$, thus proving the claim.

We now appeal to the Castelnuovo-Enriques classification of algebraic
surfaces. By
Corollary~\ref{cor5.6}, the Kodaira dimension $\kappa(X) = -1$ or $0$. If
$\kappa(X)
= -1$ then $X$ is either rational or ruled \cite[Theorem V.6.1]{Hart}. By
Proposition~\ref{prop5.11}, $X$ cannot be rational and by
Lemma~\ref{lem.ruled}, $X$
cannot be ruled. Thus we may assume without loss of generality that
$\kappa(X) = 0$.
Here there are four possibilities: (1) a $K3$ surface, (2) an Enriques
surface, (3) an
abelian surface, or (4) a hyperelliptic surface; see \cite[p. 373]{BH} or
\cite[Theorem V.6.3]{Hart}. Of these four types, only (4) has $p_g(X) = 0$
and $q = 1$. On the other hand, 
by Lemma~\ref{lem.hyperelliptic}, $X$ cannot be
hyperelliptic.
This shows that $p_g(X) = 0$ and $q(X) = 1$ is impossible, and the proof
of the
theorem is complete.
\end{proof}

\section{Wild automorphisms of abelian varieties}
\label{sec6}

In this section we will classify wild automorphisms of abelian varieties.
We will assume throughout that the base field $k$ is algebraically closed
and of
characteristic zero. We first review some basic definitions and
facts from the theory of abelian varieties; see \cite{Mu}
and \cite{La1} for detailed treatments of this material.
We will use additive notation for the group law on an abelian variety.

Let $X$ and $Y$ be arbitrary abelian varieties throughout the following
discussion.  We
write $\caphom(X, Y)$ for the group of
\emph{homomorphisms} from $X$ to $Y$; homomorphisms are,
by definition, the regular morphisms preserving the group structure.  
We also write $\End(X)$ for 
$\caphom(X,X)$, the ring of \emph{endomorphisms} of $X$. We will
use the words
\emph{automorphism} and \emph{morphism} for regular maps
which are not necessarily assumed to be homomorphisms of groups; for
example, for any $a \in X$ one has the translation automorphism $T_a:
x \mapsto x+a$.  In fact, every
regular map $\sigma: X\to Y$ between abelian varieties is of the form
$\sigma = T_b
\cdot \alpha$, where $\alpha \in \Hom(X, Y)$ and $b \in Y$; see
\cite[Theorem 4, p. 24]{La1}.
We say that $\psi \in \caphom(X, Y)$ is an \emph{isogeny}
if $\psi$ is surjective and has finite kernel. For example, for any
abelian variety $X$ and $p \in \mb{Z}$ the map $p \cdot Id_X: X \to X$
defined by $x \mapsto px$ is an isogeny \cite[Corollary IV.2.1]{La1}.
If there exists an isogeny from $X$ to $Y$, then there also must exist
an isogeny from $Y$ to $X$ \cite[p. 29]{La1}, and $X$ and $Y$ are said
to be \emph{isogeneous}.

A fundamental result is the \emph{complete reducibility} theorem of
Poincar{\'e}: if $Z$ is an abelian subvariety of $X$, then there exists
another abelian subvariety $Z' \subseteq X$ such that $Z + Z' = X$ and
$Z \cap Z'$ is a finite group \cite[Theorem II.1.6]{La1}; in this case
clearly $X$ is isogeneous to $Z \times Z'$. An abelian variety $X$ is
\emph{simple} if $0$ and $X$ are the only closed irreducible subgroups
of $X$. Abelian varieties satisfy unique decomposition into simples
in the following sense:  every abelian variety $X$ is isogeneous to
a product $X_1 \times X_2 \times \dots \times X_n$ where each $X_i$ is
a simple abelian variety, and the simples appearing in such a
decomposition are unique up to isogeny and order of the factors
\cite[p. 30]{La1}.
If $X'$ is a closed subgroup of an abelian variety $X$,
then the factor group $X/X'$ has a natural structure of an abelian
variety \cite[p. 3]{La1}.

Next, we recall some results about the structure of the endomorphism ring
$\End(X)$.
Let $X$ be isogeneous to the product $\prod_{i=1}^h X^{n_i}_i$ where the
$X_i$ are
simple abelian varieties, mutually non-isogeneous. By \cite[Corollary
VII.1.2]{La1},
$\End(X)$ is a torsionfree $\mb{Z}$-module, and so there is an embedding
of rings
$\End(X)\subset \End_{\mb{Q}}(X) := \End(X) \otimes_{\mathbb{Z}} \mb{Q}$.
Then by \cite[Theorem 7, p. 30]{La1},  each $D_i = \End_{\mb{Q}}(X_i)$ is a
finite-dimensional division algebra over $\mb{Q}$, and
$$
\End_{\mb{Q}}(X) \cong \prod_{i=1}^h M_{n_i}(D_i).
$$
In particular, if $X$ and $Y$ are isogeneous then $\End_{\mb{Q}}(X)
\cong \End_{\mb{Q}}(Y)$. If $R$ is a ring with identity, we call an
element $r \in R$
{\it unipotent} if $r-1$ is nilpotent and {\it quasi-unipotent} if some
power of $r$
is unipotent.  We shall apply these terms in particular to endomorphisms
$\alpha \in
\End(X)$.   (The reader should not confuse these definitions with the use
of the same
terms in the theory of twisted homogeneous coordinate rings in \cite{Ke};
see the
discussion in \S\ref{sec7} below.) Finally,

\begin{definition} \label{def.generate}
The algebraic subgroup of an abelian variety $X$
\emph{generated} by $S \subset X$ is the Zariski closure of
the (abstract)
subgroup $\mathopen< S \mathclose>$ of $X$ generated by $S$.
In particular, we say that $S$ {\em generates} $X$ if $S$ is not contained
in any proper closed subgroup of $X$.  If $S = \{ s_1, s_2, \dots \}$,
then
we will sometimes say that $s_1, s_2, \dots$ generate $X$.
\end{definition}

We are now ready to prove Theorem~\ref{thm0.2}(a) from the introduction,
which gives a characterization of wild automorphisms of abelian varieties.
We restate it here in a somewhat expanded form.

\begin{theorem}
\label{thm6.4}
Let $\sigma= T_b \cdot \alpha$ be an automorphism of $X$, where $\alpha
\in \End(X)$ is an automorphism and $b \in X$. Let $\beta=\alpha-Id$, and
set $S =  \{b, \beta(b), \beta^2(b),\dots \} \subset X$.
Then the following are equivalent:
\begin{enumerate}
\item
$\sigma$ is wild.
\item $\alpha$ is unipotent and $S$ generates $X$.
\item $\alpha$ is unipotent and the image $\overline{b}$ of $b$ generates
$X/\beta(X)$.
\end{enumerate}
\end{theorem}

\begin{proof} Let $Z$ be the algebraic subgroup of $X$ generated by $S$.

\smallskip
(b) $\Leftrightarrow$ (c).
If $S$ generates $X$ then clearly $\overline{b}$ generates $X/\beta(X)$.
For the converse,
assume that $\overline{b}$ generates $X/\beta(X)$.  Then $Z + \beta(X) =
X$.
Applying $\beta$ on both sides, we obtain $\beta (Z) + \beta^2(X) =
\beta(X)$. Now
\[ X = Z + \beta(X) = Z  + \beta(Z) + \beta^2(X) =  Z + \beta^2(X) \, ;\]
the last equality follows from the fact that $\beta(Z) \subseteq Z$.
Continuing in this manner, we obtain
\[ \text{$X = Z + \beta^i(X)$ for $i = 1, 2, \dots$}   \]
Since $\beta$ is nilpotent, we conclude that $Z = X$, as desired.

\smallskip
(a) $\Rightarrow$ (b):
Suppose that $\alpha-Id$ is not nilpotent. Let $Y_s$ be the connected
component containing $0$ of the closed subgroup $\ker (\alpha-Id)^s$ of
$X$.
Then $Y_s \subset Y_{s+1}$ for all $s$,
and this sequence of irreducible closed sets must stabilize, so we may set
$Y = Y_s$ for $s \gg 0$.  Clearly $\alpha(Y)\subset Y$. Hence $\alpha$
induces an endomorphism of the abelian variety $X/Y$, denoted by
$\overline{\alpha}$. By the definition of $Y$, $\overline{\alpha}-Id$ is
an
isogeny of $X/Y$. Let $\pi$ be the quotient map $X\to X/Y$ and let
$\overline{\sigma} =T_{\pi(b)} \cdot \overline{\alpha}$, which is an
automorphism of $X/Y$. Since $\overline{\alpha}-Id$ is an isogeny, it is
surjective,
and so there exists some $\overline{x}$ such that
$(\overline{\alpha}-Id)(\overline{x}) = - \pi(b)$.
Then $\overline{x} \in X/Y$ is a fixed point of the automorphism
$\overline{\sigma}$.
Set $W = \pi^{-1}(\overline{x}) \subset X$. Then $\sigma(W) = W$, and
$W \neq X$
since $Y \neq X$.
Hence $\sigma$ is not wild, a contradiction.

Thus we have proved that $\alpha$ is unipotent.  Since
$\sigma(Z)\subset Z$, and $\sigma$ is wild, we conclude also that $Z = X$.

\smallskip
(b, c) $\Rightarrow$ (a): We claim that for the purpose of
this proof we may replace $\sigma$ by $\sigma^n$, where $n \ge 1$.
Indeed, on the one hand, $\sigma$ is wild if and only if $\sigma^n$ 
is wild; see Lemma~\ref{lem2.3}(d).
On the other hand, set $\beta = \alpha - Id$ and 
\[ \gamma = \sum_{i=0}^{n-1}\alpha^i = 
\sum_{i=1}^n \binom{n}{i} \beta^{i-1} \, . \]
Then $\sigma^n=T_{b'} \cdot \alpha^n$, where $b' = \gamma(b)$.
Since $\beta$ is nilpotent, $\gamma \colon X \lra X$ is surjective 
and $\beta' =  \alpha^n-Id = \beta \gamma$ is nilpotent.
Consequently, $\beta'(X) = \beta(X)$ and if $\overline{b}$ is the 
the image of $b$ in $X/\beta(X) = X/\beta'(X)$, then
the image of $b'$ in $X/\beta(X) = X/\beta'(X)$ is
$\overline{b'} = n \overline{b}$. Consequently, 
$\overline{b'}$ generates $X/\beta'(X)$ if and only if
$\overline{b}$ generates $X/\beta(X)$. This proves the claim.

Let $W \subset X$ be a nonempty subvariety of minimal possible dimension
such that $\sigma(W) \subset W$. Our goal is to show that $W = X$;
this will imply that $\sigma$ is wild. Replacing $\sigma$ by $\sigma^n$
if necessary, we may assume also that $W$ is irreducible. Note that by
the minimality of $\dim(W)$, $\sigma$ restricts to a wild automorphism
of $W$.  By Corollary~\ref{wild auts of subvarieties}, $W$
is a translate $H + x = T_x(H)$ of some abelian
subvariety $H$ of $X$.  Now set
$$\sigma'=T_{-x}\sigma T_x = T_{b + \beta(x)} \cdot \alpha.$$
Then $\sigma'(H)=H$, so that $c = \sigma'(0) = b + \beta(x) \in H$.  Thus
$\alpha(H) \subset H$, from which it follows that $\beta(H) \subset H$.
In particular, $H$ contains $S' = \{ c, \beta(c), \beta^2(c), \dots \}$.
Since
the image of $c$ in $X/\beta(X)$ coincides with $\overline{b}$, the
implication (c)
$\Rightarrow$ (b) proved above shows that $S'$ generates $X$. Since $S'
\subset H$, we conclude that $H=X$ and hence $W = X$, as we wished to
show.
\end{proof}

For the remainder of this section we will study the conditions
for wildness given by Theorem~\ref{thm6.4} a bit further.
In view of part (c) of that theorem, we would
like to understand when a single element of an abelian
variety will generate that variety.  We analyze this situation
next.  The first two parts of the proposition
below show that if the base field $k$ is uncountable, then
``most'' points $a \in X$ will
generate $X$.  Since an arbitrary
$X$ is isogeneous to some product of simple abelian varieties, the
remaining parts
may be combined to give more precise information about which points $a \in
X$ generate $X$.
\begin{proposition}
\label{prop6.2}
Let $X$ be an abelian variety and $a \in X$ some point.
\begin{enumerate}
\item The element $a$ generates $X$ if and only if $f(a) \neq 0$ for every
$0 \neq f \in \End(X)$.
\item There is a countable set of closed subgroups $\{ G_{\alpha} \}$ of
$X$
such that $a$ generates $X$ if and only if $a \not \in \bigcup
G_{\alpha}$.
\item
If $f: X\to Y$ is an isogeny, then $a$ generates $X$
if and only if $f(a)$ generates $Y$.
\item
If $X$ is simple, then $a$ generates $X$ if and only if
$a$ is a point of infinite order on $X$.
\item
Let $X=X_1 \times X_2 \times \dots \times X_n$ where the $X_i$ are
abelian varieties such that $\Hom(X_i,X_j)=0$ for all $i\neq j$.
Then $a = (a_1, a_2, \dots, a_n)$ generates $X$
if and only if $a_i$ generates $X_i$ for every $i = 1, \dots, n$.
\item
$a=(a_1,\dots, a_n)\in X^{\times n}$
generates $X^{\times n}$ if and only if the following condition holds:
given endomorphisms $\theta_i\in \End(X)$ with $\sum_{i=1}^n
\theta_i(a_i)=0$, one must have
$\theta_i=0$ for $i=1,\dots, n$.
\end{enumerate}
\end{proposition}
\begin{proof}
\smallskip
(a) If $f(a) = 0$ for some $0\neq f \in \End(X)$
then the algebraic subgroup $Z$ generated by $a$ is contained in
$\operatorname{Ker}(f) \neq X$, so that $a$ does not generate $X$.

Conversely, suppose that $a$ does not generate $X$, so that
$a$ is contained in a closed subgroup $Z \subsetneq X$.
Let $Z_0$ be the connected component of $Z$ containing the identity
element $0$; then $pa \in Z_0$
for some integer $p \ge 1$.
By complete reducibility, there is a non-trivial abelian subvariety
$W \subset X$ such that $X= Z_0 + W$ and $Z_0 \cap W$ is finite.  Then
we may choose some isogeny $\theta: X/Z_0 \to W$, and thus we have
a nonzero endomorphism
$$
f:\quad  X \overset{p \cdot Id_X}{\lra} X \lra X/Z_0
\overset{\theta}{\lra} W \hookrightarrow X
$$
such that $f(a) = 0$.

(b) For every $0 \neq \alpha \in \End(X)$, let $G_\alpha =\ker \alpha$. By
\cite[Theorem IV.19.3]{Mu}, $\End(X)$ is countable, so part (a) implies
the result.

(c) Suppose that $a$ generates a closed subgroup $A$ in $X$ and that
$f(a)$ generates
$B$ in $Y$. Then clearly $f(A) = B$. Now it is easy to see that
$$A = X\Longleftrightarrow \dim(A) = \dim(X)
\Longleftrightarrow \dim(B) = \dim(Y) \Longleftrightarrow B = Y.$$

(d) Suppose $a$ does not generate $X$. Then, by part (a),
$a\in \operatorname{Ker}(f)$, for some $0 \ne f \in \End(X)$. 
Since $X$ is simple,
$f$ is an isogeny. Hence, $\operatorname{Ker}(f)$ is finite, 
and $a$ has finite order.
Conversely, if $a$ is a point on $X$ of finite
order then $a$ generates a finite subgroup $F$ of $X$ and
clearly, $F \ne X$.

(e) This is an easy consequence of the fact that every endomorphism
of $X$ maps $X_i$ to $X_i$ for all $i$.

(f) Suppose that $a$ does not generate $X^{\times n}$. Then by part (a),
there exists an $0 \neq \alpha\in \End(X^{\times n})$ such that $\alpha(a)=0$.
If $Z \subset \im \alpha$ is a simple abelian subvariety, then $Z$
is isogeneous to some simple abelian variety appearing in the
decomposition
of $X$ into simples.  It follows that $\Hom(\im \alpha, X) \neq 0$.
Then picking a nonzero homomorphism $f: \im \alpha \to  X$, we see that
$\theta:=f\alpha: X^{\times n} \to X$ is a nonzero homomorphism such that
$\theta(a)=0$.
But $\theta(a)=\sum_{i=1}^n \theta_i(a_i)$ for some $\theta_i\in \End(X)$,
at
least one of which is nonzero.
The converse is proved by reversing this argument.
\end{proof}

Let us also make some remarks concerning unipotent automorphisms of
abelian
varieties.  Of course, any abelian variety $X$ has at least the unipotent
automorphism $\sigma = Id$. In the next result we will identify those $X$
for which
there exist non-identity unipotent automorphisms, and show how to
construct some of
them.  Suppose that $X = Y^{\times n}$ where $Y$ is an abelian variety.
Then for any
integer matrix $M \in M_n(\mb{Z})$ we may define an endomorphism $\alpha_M
\in
\End(X)$, as follows.  Write an arbitrary point in $X$ as a column vector
$\overline{x}$ with entries from $Y$.  Then let $\alpha_M$ be defined by
the formula
$\alpha_M(\overline{x}) = M \overline{x}$, where the right hand side is
``matrix
multiplication'' performed using the $\mb{Z}$-module structure of $Y$.
\begin{proposition}
\label{unipotent auts}
\begin{enumerate}
\item Let $X = Y^{\times n}$ for some abelian variety $Y$.  Then for $M
\in
\GL_n(\mb{Z})$, the automorphism
$\alpha_M \in \End(X)$ is (quasi)-unipotent if and only if $M$ is a
(quasi)-unipotent
matrix in $\GL_n(\mb{Z})$.
\item Let $X$ be an abelian variety which is isogeneous to a product
$\prod
_i X_i$, where the $X_i$ are
simple abelian varieties.  Then $X$ has a unipotent automorphism
$Id \neq \alpha \in \End(X)$ if and only if $X_i$ and $X_j$ are isogeneous
for some $i \neq j$.
\end{enumerate}
\end{proposition}
\begin{proof}
(a) Since the mapping $\psi: M \to \alpha_M$ is a ring homomorphism from
$M_n(\mb{Z})$ to $\End(X)$, it is clear
that $M$ unipotent implies $\alpha_M$ unipotent.  Conversely, suppose that
$M$ is
not unipotent, so $M - Id$ is not
nilpotent. One may easily check that $\ker \psi = 0$, and thus
$\alpha_{M -Id} = \alpha_M - Id$ is not nilpotent in $\End(X)$, so
$\alpha_M$ is not unipotent.  Thus $M$ is unipotent if and only if
$\alpha_M$ is
unipotent; the same statement for quasi-unipotence follows immediately.

(b) Suppose first that the $X_i$ are pairwise non-isogeneous.
Then $\End(X) \subset \End_{\mb{Q}}(X)$, where $\End_{\mb{Q}}(X) \cong
\prod \End_{\mb{Q}}(X_i)$ is a
product of division rings, and hence every nilpotent element in this ring
is zero.
Now if $\alpha \in \End(X)$ is a unipotent automorphism, then $\alpha - Id
= 0$ and so $\alpha = Id$.

Conversely, if $X_i$ and $X_j$ are isogeneous for some $i \neq j$, then
the ring
$R := \End_{\mb{Q}}(X_i \times X_i) \cong \End_{\mb{Q}}(X_i \times X_j)$
embeds in $\End_{\mb{Q}}(X)$.
Here $R \cong M_2(D)$, where $D = \End_{\mb{Q}}(X_i)$ is a division ring
finite over $\mb{Q}$.
By part (a), $R$ contains many unipotent automorphisms of
the form $\alpha_M$ for non-identity unipotent matrices $M \in
\GL_2(\mb{Z})$.
\end{proof}

\section{Unipotency on the Neron-Severi Group}
\label{sec7}
As usual, we assume in this section that our base field $k$ is
algebraically closed of
characteristic zero.

In the last section we saw how to construct wild automorphisms $\sigma$ of
an abelian variety $X$.  Then by Proposition~\ref{prop2.2}, we see that we
obtain lots of examples of projectively simple twisted homogeneous
coordinate
rings $B(X, \mc{L}, \sigma)$, as long as we can find $\sigma$-ample
invertible sheaves $\mc{L}$
on $X$.  We shall show in this section that it will suffice to take
$\mc{L}$ to be any ample
invertible sheaf on $X$.  To do this, we apply criteria of Keeler for
$\sigma$-ampleness, which will
require us to consider the induced action of automorphisms
of $X$ on the group of divisors modulo numerical equivalence.

Let us begin with a review of Keeler's results from \cite{Ke}.
Let $X$ be any projective scheme, and let $\Pic X$
be the Picard group of all invertible sheaves on $X$. Two invertible
sheaves $\mc{L}$ and $\mc{L}'$ on $X$ are \emph{numerically equivalent}
if $(\mc{L} . C) = (\mc{L}' . C)$ for all integral curves $C \subset X$.
Define the
group $\Num(X)$ to be $\Pic X$ modulo numerical equivalence; it is
isomorphic to $\mb{Z}^m$ for some $m$.
When $X$ is an abelian variety, then $\Num(X)$ is isomorphic to
the Neron-Severi group $NS(X)$ \cite[p. 367]{Hart}.
Any morphism $\sigma: X \to X$ naturally induces via
pullback an endomorphism of the group $\Num(X)$, or equivalently a
matrix $P_{\sigma} \in M_m(\mb{Z})$.

Recall that we defined the notions of unipotence and quasi-unipotence
for elements of a ring in the last section.  Now let $X$ be an
projective scheme and $\sigma: X \to X$ any
automorphism.  Then we will call $\sigma$ \emph{Num-unipotent} if the
corresponding matrix $P_{\sigma} \in M_m(\mb{Z})$ is unipotent and
\emph{Num-quasi-unipotent} if $P_{\sigma} \in M_m(\mb{Z})$ is
quasi-unipotent.
(Note that in \cite{Ke} $\sigma$ is called simply unipotent or
quasi-unipotent,
respectively.)  Keeler proved that if $\sigma$ is Num-quasi-unipotent then
every ample invertible sheaf $\mc{L}$ on $X$ is $\sigma$-ample, and if
$\sigma$
is not Num-quasi-unipotent then no $\sigma$-ample invertible sheaf exists
\cite[Theorem 1.2]{Ke}.  Thus to construct $\sigma$-ample sheaves on
abelian varieties,
we just need to understand which automorphisms of abelian varieties are
Num-quasi-unipotent.

Let $X$ be an abelian variety.
For every invertible sheaf $\mc{L} \in
\Pic X$
and every $a \in X$, let $\phi_{\mc{L}}(a) = (T_a)^*(\mc{L}) \otimes
\mc{L}^{-1} \in \Pic X$.
The \emph{Theorem of the Square} \cite[Corollary 4, p. 59]{Mu} states that
$$ (T_{a+b})^*(\mc{L}) \otimes \mc{L} \cong (T_a)^*(\mc{L}) \otimes
(T_b)^*(\mc{L})\ \quad \text{for all}\ a, b \in X.$$
It follows that $\phi_{\mc{L}}$ is a group homomorphism from $X$ to
$\Pic(X)$.
Now we can easily handle the case of translation automorphisms.
\begin{lemma}
\label{lem7.1}
Let $X$ be an abelian variety with $a \in X$ and let $T_a$ be the
corresponding translation automorphism.  Then the induced automorphism
$P_{T_a}$ of $\Num(X)$ is the identity.  In particular, $T_a$ is
Num-unipotent.
\end{lemma}
\begin{proof}
The \emph{Picard variety} of $X$, denoted by $\Pic^0(X)$, is the subgroup
of $\Pic(X)$ consisting of all $\mc{L}$ such that $(T_a)^*(\mc{L}) \cong
\mc{L}$ for all $a \in X$,
in other words those $\mc{L}$ such that $\phi_{\mc{L}}$ is identically
$0$.
Then another application of the Theorem of the Square shows that for every
$\mc{L} \in \Pic(X)$ and $a \in X$,
$\phi_{\mc{L}}(a) = (T_a)^*(\mc{L}) \otimes \mc{L}^{-1} \in \Pic^0(X)$
(see \cite[p. 74]{Mu}).  Since any divisor in $\Pic^0(X)$ is numerically
equivalent to $0$~\cite[Proposition IV.3.4]{La1},
$\mc{L}$ and $(T_a)^*(\mc{L})$ induce the same element of $\Num(X)$, and
so it follows that $P_{T_a} = Id$
as required.
\end{proof}

Since an arbitrary automorphism $\sigma \in \Aut(X)$ has the form $\sigma
= T_b \cdot \alpha$ for $b \in X$ and
$\alpha \in \End(X)$, we see from the preceding result that $\sigma$ and
$\alpha$ induce the same map on $\Num(X)$.  Thus from now on we may
restrict our attention
to automorphisms in $\End(X)$.
We would like to say that if $\alpha \in \End(X)$ is a quasi-unipotent
automorphism, then
$P_{\alpha}$ is a quasi-unipotent automorphism of $\Num(X)$.  This is not
obvious, because
the map $P: \End(X) \to \End(\Num(X))$
defined by $\alpha \mapsto P_{\alpha}$ is only a homomorphism of
multiplicative semigroups, not a ring map.  Nevertheless, with a little
more work we will
succeed in proving
that $P$ preserves quasi-unipotency in Proposition~\ref{prop7.4} below.
The main step is to
express unipotency purely in terms of the multiplicative structure, which
will require the following technical bit
of algebra.
\begin{lemma}
\label{lem7.3}
Let $D_1, D_2, \dots, D_h$ be division algebras over $\mb{Q}$.  Fix
some integers $n_i>0$ and let $M \in \prod_{i=1}^{h} M_{n_i}(D_i)$.
\begin{enumerate}
\item
If $M$ is unipotent, then $M^p$ is conjugate to $M$ for all $p > 0$.
If $M$ is quasi-unipotent, then $M^p$ is conjugate to $M^q$ for some
$0<p<q$.
\item
Assume that each $D_i$ is a (commutative) field.  Then the converses of
both statements
in part (a) hold.
\item
Let $g: \prod_{i=1}^{h} M_{n_i}(D_i)\to M_m(F)$ be a homomorphism of
multiplicative semigroups, where $F$ is a field.  Then $g$ preserves
quasi-unipotency (respectively, unipotency).
\end{enumerate}
\end{lemma}

\begin{proof} In both parts it is easy to reduce to the case where
$h=1$, and we do so, writing $D = D_1$ and $n = n_1$.

(a) Let $M \in M_n(D)$ be unipotent.  It is an exercise in linear algebra
that in this case $M$ has a Jordan canonical form (this is of course not
necessarily true for arbitrary matrices in $M_n(D)$).  Also, $M$ and $M^p$
must have the same Jordan form for all $p > 0$, so they are conjugate.

Suppose instead that $M$ is quasi-unipotent.  Then $M^s$ is unipotent for
some $s>0$, and so $M^s$ and $M^{sp}$ are conjugate for all $p > 0$.

(b) Assume now that $D$ is commutative, and let $M \in M_n(D)$.  Let
$\lambda_1, \cdots, \lambda_n$ be the eigenvalues of $M$ in some algebraic
closure of $D$.  Then $M$ is quasi-unipotent if and only if all of the
$\lambda_i$ are roots of unity.

Suppose that $M^p$ is conjugate to $M^q$ for some $0<p<q$.  Then
the ordered set $\{\lambda_1^p,\cdots, \lambda_n^p\}$ is a permutation
of the ordered set $\{\lambda_1^q,\cdots, \lambda_n^q\}$. Let $\tau
\in S_n$ be that permutation, namely, $\lambda_i^p=\lambda_{\tau(i)}^q$
for all $i$. If $w=n!$, then $\tau^w$ is the identity. Hence we have
$$\lambda_i^{p^w}=\lambda_{\tau(i)}^{p^{w-1}q}=\cdots=
\lambda_{\tau^t(i)}^{p^{w-t}q^t}=\cdots=\lambda_{\tau^w(i)}^{q^w}
=\lambda_i^{q^w}.$$
Thus $\lambda_i^{q^w-p^w}=1$ for all $i$ and $M$ is quasi-unipotent.

If $M$ is conjugate to $M^p$ for all $p>0$, then $M$ is quasi-unipotent
by the last paragraph. Hence $M^p$ is unipotent for some $p$. Since $M$
is conjugate to $M^p$, $M$ is unipotent.

(c) This is an immediate consequence of parts (a) and (b).
\end{proof}

\begin{remark}
\label{rem7.5}
Lemma \ref{lem7.3}(a) also allows us to improve Theorem \ref{thm5.3}
and Corollary \ref{cor5.4} slightly. Suppose that $X$ is an abelian
variety, that $\sigma \in \Aut(X)$ is wild, and that
some power of $\sigma$ is a translation. Then $\sigma = T_b \cdot
\alpha$ with $\alpha \in \End(X)$ and $\alpha^n = Id$ for some $n
\geq 0$.  By Theorem~\ref{thm6.4}, $\alpha$ is unipotent.
Then by Lemma~\ref{lem7.3}(a), $\alpha$ is conjugate to $\alpha^n = Id$,
so $\alpha = Id$. Thus the conclusion of Theorem \ref{thm5.3} and
Corollary \ref{cor5.4} may be improved to the statement that $\sigma$
is a translation.
\end{remark}

Now we return to our abelian variety $X$,
and show how the above lemma may be applied to the map $P: \End(X) \to
\End(\Num(X))$.
First, we want to know the action of $P$ on the ``scalars''.
\begin{lemma}
\label{lem7.2}
Let $X$ be an abelian variety. Let $\alpha= n \cdot Id_X \in \End(X)$
where $n \in \mb{Z}$. Then $P_\alpha = n^2 \cdot Id_{\Num(X)}$.
\end{lemma}
\begin{proof}
This is immediate from \cite[Proposition 2, p. 92]{La1}.
\end{proof}
Let $\End_{\mb{Q}}(X) = \End(X) \otimes_{\mb{Z}} \mb{Q}$ and
$\Num_{\mathbb{Q}}(X)=\Num(X)\otimes_{\mathbb{Z}}\mathbb{Q}$.  By
Lemma \ref{lem7.2}, if we define $P_{n^{-1}Id_X}=n^{-2}Id_{\Num(X)}$ for
$n \geq 1$  then $P$ naturally extends to a semigroup homomorphism
$P: \End_{\mathbb{Q}}(X)\to \End_{\mb{Q}}(\Num(X))$.

Let $X$ be isogeneous to a product $\prod_{i=1}^h
X^{n_i}_i$ where the $X_i$ are simple abelian varieties, mutually
non-isogeneous.  As we stated in \S\ref{sec6}, each $D_i =
\End_{\mathbb{Q}}(X_i)$ is a finite dimensional division algebra over
$\mb{Q}$, and $\End_{\mathbb{Q}}(X) \cong \prod_{i=1}^h M_{n_i}(D_i)$.
Also, letting $m$ be the rank of $\Num(X)$ we have
$\End_{\mb{Q}}(\Num(X)) \cong M_{m}(\mb{Q})$.
We may now prove our main result on Num-quasi-unipotence for automorphisms
of abelian varieties.  This settles Theorem~\ref{thm0.2}(b) from the
introduction.

\begin{theorem}
\label{prop7.4}
Let $X$ be an abelian variety.
\begin{enumerate}
\item
The map $P:\End_{\mathbb{Q}}(X) \to \End_{\mb{Q}}(\Num(X))$ preserves
(quasi)-unipotency.
\item
If $\alpha \in \Aut(X)$ is (quasi)-unipotent and $b \in X$ is any point,
then $\sigma = T_b \cdot \alpha \in \Aut(X)$ is Num-(quasi)-unipotent.
\item
If $\sigma$ is a wild automorphism of $X$, then
$\sigma$ is Num-unipotent and any ample invertible sheaf $\mc{L}$ on $X$
is
$\sigma$-ample.
\end{enumerate}
\end{theorem}

\begin{proof}
(a) As we saw in the comments before the theorem, the
map $P$ is a semigroup homomorphism from
$\prod_{i=1}^h M_{n_i}(D_i)$ to $M_m(\mb{Q})$.  Thus
the result is immediate from Lemma~\ref{lem7.3}(c).

(b) Since $P$ is a semigroup homomorphism, this follows from part (a) and
Lemma~\ref{lem7.1}.

(c) This follows from part (b), Theorem~\ref{thm6.4}(b), and
\cite[Theorem 1.2]{Ke}.
\end{proof}

We have finally shown all of the necessary pieces to demonstrate the
existence of a large supply of
projectively simple twisted homogeneous coordinate rings.
\begin{corollary}
Let $X$ be an abelian variety with automorphism $\sigma = T_b \cdot
\alpha$, where
$\alpha \in \End(X)$ is unipotent and the image of $b$ generates
$X/(\alpha- Id)(X)$.
Then if $\mc{L}$ is any ample invertible sheaf on $X$, then the ring $B(X,
\mc{L}, \sigma)$
is projectively simple.
\end{corollary}
\begin{proof}
By Theorem~\ref{thm6.4}, $\sigma$ is a wild automorphism.  By
Theorem~\ref{prop7.4}, the sheaf
$\mc{L}$ is $\sigma$-ample.  Then by Proposition~\ref{prop2.2},
$B(X,\mc{L}, \sigma)$ is
a projectively simple ring.
\end{proof}

In \cite[Example 5.18]{AV}, Artin and Van den Bergh gave an example of
automorphism $\sigma$ of a surface $Y$, defined by ``translation
along a pencil of elliptic curves'', which has the property that
$P_{\sigma}$ is unipotent, but not the identity.
There seem to be few
other explicit examples in the literature of automorphisms of projective
schemes which are Num-quasi-unipotent, but whose action on $\Num(X)$ is
non-trivial.  In part
(c) of the next result we show
that we can produce numerous such examples in the setting of abelian
varieties.

\begin{lemma}
\label{lem7.6}
Let $X$ be an abelian variety.
\begin{enumerate}
\item
If $f\in \End(X)$ is an isogeny, then $P_f$ is nonzero.
\item
If $\beta\in \End_{\mathbb{Q}}(X)$ is nonzero, then $P_\beta$ is nonzero.
\item
If $\alpha\in \End(X)$ is unipotent and $\alpha \neq Id$, then $P_\alpha$
is
unipotent and $P_{\alpha} \neq Id$.
\end{enumerate}
\end{lemma}

\begin{proof}
(a) This follows from Lemma \ref{lem7.2} and the fact that there is a
$g\in \End(X)$ such that $fg =n \cdot Id_X$ \cite[p. 29]{La1}.

(b) Replacing $\beta$ by $n\beta$ for some $n \geq 1$ we may assume that
$\beta\in \End(X)$. If $Y$ is the image of $\beta$, then there are
endomorphisms $f: Y\to X$ and $g: X\to Y$ such that $g \beta f\in \End(Y)$
is an isogeny. Now the assertion follows from (a).

(c) By Theorem~\ref{prop7.4}, $P_\alpha$ is unipotent, so we just need to
show that $P_\alpha$ is not the identity. Let $\beta=\alpha-Id$, so that
$\beta$ is nilpotent. Let $p$ be the maximal integer such that $\beta^p
\neq 0$. Then
$$\alpha^n=\sum_{i= 0}^p \binom{n}{i} \beta^i$$
for all $n \geq 0$. By \cite[Proposition 2, p. 92]{La1},
\begin{equation}
\label{E1}
P_{\alpha^n} =\frac{1}{2}\sum_{i= 0}^p \sum_{j = 0}^p
\binom{n}{i}\binom{n}{j}
D(\beta^i, \beta^j)
\end{equation}
for all $n \geq 0$, where $D(f,g)=P_{f+g}-P_{f}-P_{g}$ for all
$f,g\in \End(X)$.  For all $n$
the right hand side of \eqref{E1} is equal to $Q(n)$, where $Q$
is some fixed polynomial
$$
Q(z) = \sum_{i = 0}^{2p} h_i z^i \in \End(\Num(X))[z] = M_m(\mb{Z})[z].
$$
Now suppose that $P_{\alpha}$ is the identity.  Then $P_{\alpha^n}$ is the
identity for
all $n \geq 1$, and thus $Q(n) = Id$ for all $n \geq 1$.  Let $| M |$ be
notation for the norm of an arbitrary matrix $M \in M_{m}(\mb{Z})$ in the
Euclidean topology
on $\mb{R}^{m^2}$.
Now if $h_i \neq 0$ for some $i \geq 1$, then clearly $\lim_{n \to \infty}
| Q(n)| = \infty$,
contradicting the fact that $Q(n) = Id$ for all $n \geq 1$.  Thus $h_i =
0$ for all $i \geq 1$.

In particular, the leading coefficient $h_{2p}$ of $Q$ is
$(p!)^{-2}D(\beta^p,
\beta^p) = 0$.  But using Lemma \ref{lem7.2},
$$
D(\beta^p, \beta^p) = P_{2\beta^p} - 2 P_{\beta^p} =
4P_{\beta^p}-2P_{\beta^p}
= 2P_{\beta^p},
$$
so $P_{\beta^p} = 0$ (recall our standing assumption that $\ch k = 0$).
This
contradicts part (b).  Thus $P_{\alpha} \neq Id$.
\end{proof}

In order to illustrate the results of this section, we conclude with a
simple example where the action of automorphisms on $\Num(X)$ may be
calculated explicitly.
\begin{proposition}
\label{E times E}
Let $E$ be an elliptic curve without complex multiplication (in other
words $\End(E) \cong \mb{Z}$), and let $X=E\times E$.
Let $M = \begin{pmatrix} a&b\\ c&d\end{pmatrix}\in \GL_2(\mathbb{{Z}})$ be
arbitrary and let $\alpha = \alpha_M \in \End(X)$ be the corresponding
automorphism; explicitly, $\alpha(x, y) = (ax + by, cx+dy)$.
\begin{enumerate}
\item  The divisors $C^1= 0 \times E$, $C^2= E \times 0$,
and $C^3= \{(x, x) \mid x \in E \}$ form
a basis for $\Num(X) \otimes_{\mb{Z}} \mb{Q} \cong \mb{Q}^3$.
\item
With respect to the basis given in (a), $P_\alpha \in M_3(\mathbb{{Z}})$
is
equal to
$$\begin{pmatrix} a^2-ab&  c^2-cd& (a+c)^2-(a+c)(b+d) \\
b^2-ab& d^2-cd& (b+d)^2-(a+c)(b+d)\\
 ab & cd & (a+c)(b+d)
\end{pmatrix}.$$
\item
$\det P_{\alpha}=(\det M)^3$.
\item
$\alpha$ is quasi-unipotent if and only if
$P_\alpha$ is.
\end{enumerate}
\end{proposition}
\begin{proof}
(a, b) These parts are the statements of \cite[Exercise II.4.16.2]{kollar}
and \cite[Exercise II.4.16.6]{kollar}),
respectively, and we omit the proofs.

(c) This is an easy computation using part (b).

(d) By part (b), all entries of $P_\alpha$ are polynomial functions of
$a,b,c,d$. Hence the semigroup homomorphism $P: GL_2(\mathbb{Z})\to
GL_3(\mathbb{Z})$ can be extended
to a group homomorphism $P: GL_2(\mathbb{C})\to GL_3(\mathbb{C})$
with the same formula given in (b). Since $P$ is a group homomorphism,
the desired assertion holds if and only if it holds after conjugation.

By conjugation we may assume that $\alpha$ has one of the two forms
$\begin{pmatrix} a&0
\\0 &d\end{pmatrix}$ or $\begin{pmatrix} a&1\\ 0&a\end{pmatrix}$. Then one
may check using the formula in part (b) that
the eigenvalues of $P_\alpha$ are either $\{a^2, d^2, ad\}$ or $\{a^2,
a^2, a^2\}$, respectively.
Hence $P_\alpha$ is quasi-unipotent
if and only if $M$ is a quasi-unipotent matrix, if and only if $\alpha_M$
is
a quasi-unipotent, by Proposition~\ref{unipotent auts}(a).
\end{proof}

\section{Classification of low-dimensional cases}
\label{sec8}

In this section we tie together our previous results to classify
projectively simple twisted homogeneous coordinate rings of small
GK-dimension.  Again we
assume that $k$ is an algebraically closed field of characteristic zero.

We will use the following result of Keeler about the GK-dimension of
twisted homogeneous coordinate rings \cite[Theorem 6.1]{Ke}.  Let
$B = B(X, \mc{L}, \sigma)$ where $\mc{L}$ is $\sigma$-ample.  Then $\GKdim
B$ is an integer satisfying the inequalities
\begin{equation}
\label{GK B}
j + \dim X + 1 \leq \GKdim B \leq j(\dim X -1) + \dim X + 1,
\end{equation}
where $j + 1$ is the size of the largest Jordan block of
the matrix $P_{\sigma} \in \End(\Num(X))$ giving the induced action of
$\sigma$ on $\Num X$.  We note that if $\dim X=1$, then $\Num(X) = \mb{Z}$ and 
so $j=0$, while if $\dim X=2$, then $j=0$ or $2$ by \cite[Lemma 5.4]{AV}.  Moreover, 
it is a fact that $j$ is always even \cite[Lemma 6.12]{Ke}.

\begin{proposition}
\label{prop8.1}
Let $B = B(X, \mc{L}, \sigma)$, where $\mc{L}$ is $\sigma$-ample.
Then $B$ is a projectively simple ring with $\GKdim(B) = 2$ or $3$
if and only if
\begin{enumerate}
\item
$X$ is an abelian variety of dimension $\GKdim(B) - 1$,
\item
$\sigma$ is the the translation $T_b$, and
\item
$b$ generates $X$.
\end{enumerate}
\end{proposition}

\begin{proof}
Suppose first that $B$ is projectively simple.
Since $\GKdim B \leq 3$, by \eqref{GK B} we must have $\dim X \leq 2$, and
since
$j$ in that equation is even, the only possibility is $j = 0$, which
forces
$\GKdim(B) = \dim X + 1$.  Also, $\sigma$ is wild by
Proposition~\ref{prop2.2}. By Theorem~\ref{thm5.12},
$X$ must be an abelian variety.  Let $\sigma=T_b \cdot \alpha$ for
$\alpha \in \End(X)$ and $b \in X$.  Note that $P_{\sigma} = P_{\alpha}$
by
Lemma \ref{lem7.1}. By Theorem~\ref{thm6.4}, $\alpha$ is
unipotent.  If $\alpha$ is not the identity, then by
Lemma~\ref{lem7.6}(c),
$P_\sigma = P_\alpha$ is unipotent and not the identity.  In this case
$P_\sigma$ has a non-trivial Jordan block, and so $j \neq 0$ in
\eqref{GK B}, a contradiction.  We conclude that $\alpha=Id$ and
$\sigma=T_b$ is a translation.  By Theorem~\ref{thm6.4} again, $b$
generates
$X$.

Conversely, suppose that $X$ is an abelian variety with $\dim X = 1$ or
$2$ and that $\sigma = T_b$ where $b$ generates $X$.  By
Theorem~\ref{thm6.4},
$\sigma$ is wild.  Then $B$ is projectively simple by
Proposition~\ref{prop2.2}.
Moreover, $P_{\sigma} = Id$ by Lemma \ref{lem7.1}, so $j = 0$ in \eqref{GK
B} and
$\GKdim(B) = \dim X +1$.
\end{proof}

One can prove analogous results for other small values of dimension.
We omit the similar proof of the following.
\begin{proposition}
\label{prop8.2}
Let $B = B(X,\mathcal{L}, \sigma)$ where $X$ is an abelian variety and
$\mc{L}$ is $\sigma$-ample.
\begin{enumerate}
\item
Suppose that $\GKdim B = 4$. Then $B$ is projectively simple if and only
if
$\dim X=3$, $\sigma = T_b$, and $b$ generates $X$.
\item
Suppose that $\GKdim B=5$. Then $B$ is projectively simple if and only if
either
\begin{enumerate}
\item[(1)]
$\dim X=4$, $\sigma = T_b$, and $b$ generates $X$, or
\item[(2)]
$\dim X=2$, $\sigma=T_b, \cdot \alpha$ for some non-identity
unipotent $\alpha\in \End(X)$, and the image $\overline{b}$ of $b$
generates
$X/(\alpha-Id)(X)$.
\end{enumerate}
\end{enumerate}
\end{proposition}

\section{The Cohen-Macaulay property}
\label{xxsec10}
Our goal in the last two sections is to study some of the
homological properties of projectively simple rings and their 
associated noncommutative schemes. In particular we will prove 
Theorem~\ref{thm0.5} stated in the introduction.
In order to get interesting results, we take as our basic hypothesis in
these two sections that all algebras have a balanced dualizing complex.
This condition
holds for many important classes of algebras.

The objects of interest and the methods we employ
here are rather different from those of the preceding sections, and we
must begin with a review of the definitions of
noncommutative projective schemes \cite{AZ1}, Serre duality \cite{YZ1},
and
dualizing complexes \cite{Ye, VdB}. We assume that the reader is familiar
with the basics of derived categories and follow the notations used in
\cite{YZ1}. For example, the $n$th complex shift of a complex $X$ is
denoted by $X[n]$.

Let $A$ be a right noetherian graded $k$-algebra. Let $\Gr A$ be the
category of graded right $A$-modules and let $\gr A$ be the full
subcategory of $\Gr A$ consisting of the noetherian objects.  Recall
from \S\ref{xxsec1} that a graded $A$-module $M$ is called
\emph{torsion} if for every $m \in M$ there is some $n \geq 0$ such
that $m A_{\geq n} = 0$. Let $\Tors A$ denote the full subcategory of
$\Gr A$ consisting of all torsion objects and let $\tors A$ be the full
subcategory of noetherian torsion objects, or equivalently modules of
finite $k$-dimension. The \emph{noncommutative projective spectrum} of
$A$ is defined to be $\Proj A:=(\QGr A, {\mathcal A},s)$ where $\QGr A$
is the quotient category $\Gr A/\Tors A$, the object ${\mathcal A}$ is the
image of $A$ in $\QGr A$, and $s$ is the auto-equivalence of $\QGr A$
induced by the degree shift $M\to M(1)$ for $M\in \Gr A$.  If we want
to work with noetherian objects only, then $\proj A:=(\qgr A,
{\mathcal A},s)$ is also called the projective spectrum of
$A$.  The canonical functor from $\Gr A \to \QGr A$ (and from
$\gr A$ to $\qgr A$) is denoted by $\pi$.  If $M\in \Gr A$, we will
use the corresponding calligraphic letter ${\mathcal M}$ for $\pi(M)$
if there is no chance of confusion. For example, ${\mathcal A}=\pi(A)$.

Let $X=\Proj A$.  For ${\mathcal N}\in \QGr A$, the $i$-th cohomology
group of $\mc{N}$ is defined to be
$$
\HB^i(X, {\mathcal N})=\Ext^i_{\QGr A}({\mathcal A},{\mathcal N})
$$
for all $i\geq 0$.
Then the (right) cohomological dimension of $X$ is defined
to be
$$\cdim(X)=\max\{i \;|\; \HB^i(X,{\mathcal N})\neq 0 \text{ for some }
{\mathcal N}\in \QGr A\}.$$
Suppose that $d=\cdim(X)$ is finite.
An object $\omega\in \qgr A$ is called a {\it dualizing sheaf}
for $X$ if there is a natural isomorphism
$$\theta^0: \HB^d(X,{\mathcal N})^*\to \Hom_{\qgr A}({\mathcal N},
\omega)$$
for all ${\mathcal N}\in \qgr A$.  Here $(-)^*$ means the $k$-vector
space dual. We say that $X$ satisfies {\it Serre duality} if a dualizing
sheaf $\omega$ exists. In this case we say that $X$ is {\it classically
Cohen-Macaulay} if $\theta^0$ can be extended to a sequence of natural
isomorphisms
$$\theta^i: \HB^{d-i}(X,{\mathcal N})^*\to \Ext^i_{\qgr A}({\mathcal N},
\omega)$$
for all $\mc{N}$ and $0 \leq i \leq d$.
If $A$ is commutative, then
the concepts of dualizing sheaf and the Cohen-Macaulay property which we
have defined here agree with the usual commutative notions.

The notion of a balanced dualizing complex for a noncommutative graded
algebra was introduced by Yekutieli in \cite{Ye}. Let $\D(\Gr A)$ denote
the derived category of graded right $A$-modules. Given a complex $Y$,
we use the notation $\HS^i(Y)$ for the $i$th cohomology of $Y$. Let
$A^{\circ}$ be the opposite algebra of $A$, and let $A^{e} = A
\otimes_k A^{\circ}$. A \emph{dualizing complex} for a noetherian graded
ring $A$ is a bounded complex $R\in \D(\Gr A^{e})$ such that
\begin{enumerate}
\item $R$ has finite injective dimension over $A$ and $A^{\circ}$,
\item The $A$-bimodule $\HS^i(R)$ is noetherian on both sides for all $i
\in \mb{Z}$, and
\item The natural maps $A \to \RHom_A(R, R)$ and $A \to
\RHom_{A^{\circ}}(R, R)$
are isomorphisms in the derived category $\D(\Gr A^{e})$.
\end{enumerate}

Let $\Gamma_{\mf{m}}: \Gr A^{e} \to \Gr A^{e}$ be the \emph{torsion
functor}
$\lim_{n \to \infty} \uHom_{A}(A/\mf{m}^n, -)$.  Writing
$\mf{m}^{\circ} = (A^{\circ})_{\geq 1}$, the
functor $\Gamma_{\mf{m}^{\circ}}$
is defined similarly.
If $A$ has a dualizing complex $R$, then $R$ is called \emph{balanced} if
there are isomorphisms $\R\Gamma_{\mf{m}}(R) \cong A^{*}$ and
$\R\Gamma_{\mf{m}^{\circ}}(R) \cong A^{*}$ in $\D(\Gr A^{e})$.

Yekutieli proved that if $A$ is noetherian with a balanced dualizing
complex $R$,
then $\omega = \HS^{-(d+1)}({\mathcal R})$ is a dualizing sheaf for $\Proj
A$,
where ${\mathcal R}=\pi(R)$ \cite[4.2(4)]{YZ1}.  In addition, $X$ is
classically
Cohen-Macaulay if and only if $\omega[d+1]$ is isomorphic to ${\mathcal
R}$ in
the derived category $\D(\QGr A)$ \cite[4.2(5)]{YZ1}, or equivalently if
and only
if ${\mathcal R}$ has nonzero cohomology in only one term.

A powerful criterion of Van den Bergh \cite[6.3]{VdB} says that $A$ admits
a balanced dualizing complex if and only if $A$ satisfies the left and
right
$\chi$ condition and $X = \Proj A$ has finite left and right cohomological
dimension.  Here, $A$ is said to satisfy the (right) $\chi$ condition if
$\dim_k \uExt_A^i(A/\mf{m}, M) < \infty$ for all $M \in \gr A$ and all $i
\geq 0$, and similarly on
the left.
It follows from Van den Bergh's criterion that many important classes of
rings
admit a balanced dualizing complex, although there are examples of
noetherian
rings which do not because they fail the $\chi$ condition \cite{SZ},
\cite{KRS}.

We may now prove Theorem \ref{thm0.5}(a).
\begin{proposition}
\label{prop2.1}
Let $A$ be a connected graded, noetherian, projectively simple ring
admitting a balanced dualizing complex. Then the associated projective
scheme $\Proj A$ is classically Cohen-Macaulay.
\end{proposition}

\begin{proof} By \cite[6.3]{VdB}, $A$ satisfies the $\chi$ condition and
has finite cohomological dimension on both sides, and moreover the
balanced dualizing complex over $A$ is given by
$$R=\R \Gamma_{\mathfrak m}(A)^*.$$
Let $X = \Proj A$.  Note
that $\R^i\Gamma_{\mathfrak m}(N) \cong \HB^{i-1}(X, {\mathcal N})$
for all $N\in \Gr A$ and all $i \geq 2$ \cite[7.2(2)]{AZ1}.  Thus if
$\cd(\Proj A) = d$, then
$\R^{n}\Gamma_{\mathfrak m}=0$ for $n>d+1$, and so
$\HS^n(R)=0$ for all $n<-d-1$. Let ${\mathcal R} =\pi(R)$. Since the
dualizing sheaf $\omega=\HS^{-(d+1)} ({\mathcal R})$ is necessarily
nonzero, $M:=\HS^{-(d+1)}(R)$ is not finite dimensional over $k$.
Let $j_0=-d-1$ and put $j=\max\{i\;|\; \HS^i({\mathcal R})\neq 0\}$.
We want to show that $j_0=j$. If this is the case, then $\omega \cong
{\mathcal R}[-d-1]$ and we are done by \cite[4.2(5)]{YZ1}.

By the definition of a dualizing complex, $\HS^n(R)$ is noetherian
on both sides for every $n$.  Thus for every $n>j$, $\HS^n(R)$ is
finite dimensional over $k$ because $\HS^n({\mathcal R})=0$, but
$N:=\HS^j(R)$ is not finite dimensional since $\HS^j(\mc{R}) \neq 0$.

Let $Y$ be the truncation $\tau^{\leq j}R$ and $Z$ the truncation
$\tau^{\geq (j+1)}R$.  Then we have a distinguished triangle in the
derived category $\D(\Gr A)$
$$Y\to R\to Z\to Y[1]$$
which induces a long exact sequence
$$\cdots \to \Ext^i_A(Z,R)\to \Ext^i_A(R,R)\to \Ext^i_A(Y,R)\to \cdots.$$
Since $\HS^n(Z)$ is finite dimensional and $\HS^n(R)$ is noetherian for
all $n$, it follows from the $\chi$ condition and induction on the
lengths of the bounded complexes $Z$ and $R$ that $\Ext^i_A(Z,R)$ is
finite dimensional for all $i$.  If we show that $\Ext^i_A(Y,R)$ is
infinite dimensional for some $i$, then it will follow that
$\Ext^i_A(R,R)$ is also infinite dimensional. Since $Y$ is bounded
above at $j$ with $\HS^j(Y)=\HS^j(R)=N$ and $R$ is bounded below at
$j_0$ with $\HS^{j_0}(R)=M$, every nonzero map from $N$ to $M$ induces
a nonzero element in $\Ext^{j_0-j}_A(Y,R)$.  Thus there is an injection
of vector spaces
$$
\uHom_A(N,M) \to \Ext^{j_0-j}_A(Y, R).
$$
Now since $N$ and $M$ are not torsion, by
Lemma~\ref{lem1.7}(b) $\uHom_A(N,M)$ is infinite dimensional.  Hence
$\Ext^{j_0-j}_A(Y,R)$ is infinite dimensional, whence
$\Ext^{j_0-j}_A(R,R)$
is infinite dimensional. But since $R$ is a dualizing complex, by
definition we we have $\Ext^i_A(R,R)=0$ for all $i\neq 0$. Thus $j=j_0$.
\end{proof}

\section{The Gorenstein property}
\label{xxsec11}
In this section, we will complete the proof of Theorem~\ref{thm0.5}
by showing that if $A$ is a projectively simple ring with a balanced
dualizing
complex, then the dualizing sheaf $\omega$ must be invertible. This is
equivalent to
$X = \Proj A$ being Gorenstein in the commutative case. In preparation for
this
result, we will first study graded bimodules over two rings.

Throughout this section, let $A$ and $B$ be noetherian finitely graded
prime $k$-algebras, with graded quotient rings $Q$ and $T$ respectively.
Given two noetherian $(A, B)$-bimodules $M$ and $N$, we say $M$ and $N$
are {\it p.isomorphic} (short for projectively isomorphic), and write $M
\cong_p N$, if there exists an isomorphism of bimodules $f: M_{\geq n}
\to N_{\geq n}$ for some $n\gg 0$. The map $f$ induces an isomorphism
$\pi(M) \to \pi(N)$ in $\Qgr B$ (and similarly in $\Qgr A^{\circ}$).
A map $f:M\to N$
is called a {\it p.isomorphism} if $f_{\geq n}: M_{\geq n}\to N_{\geq n}$
is an isomorphism for some $n\gg 0$.

If $M$ is a graded right $B$-module and $N$ is a graded left $B$-module,
the graded tensor product of $M$ and $N$ over $B$ is denoted by $M
\uotimes_B N$. As usual, if $M$ is a graded $(A,B)$-bimodule then $M$ is
called \emph{invertible} if there is a graded $(B,A)$-bimodule $N$ with
bimodule isomorphisms $M \uotimes_B N \cong A$ and $N \uotimes_A M \cong
B$.
Analogously, $M$ is called {\it p.invertible} if there exists such a
graded $(B,A)$-bimodule $N$ with $M \uotimes_B N \cong_p A$ and $N
\uotimes_A M \cong_p B$, and $N$ is called the {\it p.inverse}.

Next, we want to define a different but related notion of invertibility.
We call a noetherian $(A,B)$-bimodule $M$ {\it generically invertible} if
\begin{enumerate}
\item
$M$ is Goldie torsionfree on both sides,
\item
$M$ is {\it evenly localizable} in the sense that
$$Q\uotimes_A M\cong Q\uotimes_A M\uotimes_B T \cong M\uotimes_B T, \ \
\text{and}$$
\item
$Q\uotimes_A M = M\uotimes_B T$ is an invertible $(Q,T)$-bimodule.
\end{enumerate}

Although it is useful to state all three conditions in the above
definition, we should note that in general condition (b) is actually
a consequence of condition (a), as follows.

\begin{lemma}
\label{xxlem3.1}
If $M$ is a noetherian $(A, B)$-bimodule which is Goldie torsionfree
on both sides, then $M$ is evenly localizable.
\end{lemma}

\begin{proof}
First we claim that the $(A, B)$-bimodule $M' = Q\uotimes_A M$
is Goldie torsionfree on
the right.  This is because
every finitely generated graded
right $B$-submodule of $Q\uotimes_A M$ is contained in $qA \uotimes_A
M \cong M$ for some homogeneous $q \in Q$.

Thus right multiplication by any element $b \in B$ must induce a
left $Q$-module injection $\psi_b: M' \to M'$.  Since $_Q M'$ is of
finite rank over the graded-semisimple ring $Q$, the map $\psi_b$ must
be a bijection for all $b$ and so $Q\uotimes_A M \cong Q \uotimes_A
M \uotimes_B T$. A symmetric argument gives $M \uotimes_B T  \cong
Q\uotimes_A M \uotimes_B T$.
\end{proof}

The reader may check that an invertible bimodule is generically
invertible.
The following key proposition is a kind of converse statement for
projectively simple
rings.  The proof will be given below after several more lemmas.

\begin{proposition}
\label{prop3.2}
Let $A$ and $B$ be connected graded, noetherian, projectively simple rings
admitting balanced dualizing complexes. Then every generically invertible
noetherian $(A,B)$-bimodule is p.invertible.
\end{proposition}

Let $M$ be an $(A,B)$-bimodule. Then we may define the $(B,A)$-bimodules
$^{\vee}M = \uHom_{B}(M_B,B_B)$ and $M^{\vee} =
\uHom_{A^{\circ}}(_AM,{}_AA)$.
Induced by evaluation, we get a $B$-bimodule homomorphism

\begin{equation}
\label{def.phi}
\phi_M: {}^{\vee}M \uotimes_A M\to B
\end{equation}

\noindent
and an $A$-bimodule homomorphism

\begin{equation}
\label{def.xi}
\xi_M: M\uotimes_B M^{\vee} \to A,
\end{equation}

\noindent
both of which are graded morphisms of degree $0$.
The bimodules $M^{\vee}$ and $^{\vee} M$ will be our candidates for the
p.inverse of $M$.

We shall need the following technical result, which says loosely that
``Ext commutes with localization''.  The proof is an easy
generalization to the graded setting of (\cite[Proposition 1.6]{BL}).
\begin{lemma}
\label{lem3.3} Let $B$ be a graded noetherian prime ring. Let $T$ be
the graded quotient ring of $B$. Let $M_B$ be a finitely generated
graded $B$-module and let $N$ be a graded $B$-bimodule such that
$T\uotimes_B N\cong N\uotimes_B T$. Then
$$
T \uotimes_B \uExt^i_B(M_B,{}_B N_B) \cong
\uExt^i_T(M\uotimes_B T, N \uotimes_B T)
$$
for all $i$.
\end{lemma}

We next show that to prove Proposition \ref{prop3.2} it is enough to
know that $M^{\vee}$ and $^{\vee} M$ are noetherian.

\begin{lemma}
\label{lem3.4}
Let $A$ and $B$ be two prime noetherian projectively simple rings. Let
$M$ be a generically invertible noetherian $(A,B)$-bimodule. If $M^{\vee}$
and
$^{\vee}M$ are noetherian on both sides, then $M$ is p.invertible with
p.inverse $M^{\vee}\cong_p {}^{\vee}M.$
\end{lemma}

\begin{proof} We will show that the evaluation map
$$\phi_M: {^{\vee}M}\uotimes_A M\to B$$
of \eqref{def.phi} is a p.isomorphism. By the definition of generic
invertibility, the $(Q,T)$-bimodule $L = M\uotimes_B T \cong Q\uotimes_A
M$ is invertible. By Morita theory, the inverse of $L$ must be given
by the $(T, Q)$-bimodule $L^{-1}=\uHom_T(L,T_T)$. By Lemma \ref{lem3.3},
we have
$$T\uotimes_B {^{\vee}M}\cong \uHom_T(L,T_T)=L^{-1}.$$
Hence ${^{\vee}M}$ is nonzero and the map $\phi_M$ is nonzero.
Since $B$ is projectively simple, the map $\phi_M$ must be surjective
in large degree.

Now consider the map $T \uotimes_B \phi_M$.  This map must be an
isomorphism, since it may be identified with the chain of isomorphisms
$$
T\uotimes_B {^{\vee}M}\uotimes_A M\cong L^{-1}
\uotimes_A M\cong L^{-1}\uotimes_Q (Q\uotimes_A M)\cong T.
$$
It follows that the left $B$-module $\ker \phi_M$ is Goldie torsion.
Since by assumption $^{\vee} M$ is a noetherian bimodule, $^{\vee} M
\uotimes_A M$ and thus $\ker \phi_M$ are also noetherian bimodules.
Hence $\ker \phi_M$ is finite dimensional by Lemma \ref{lem1.7}, and
$\phi_M$ is a p.isomorphism as we wished. The proof that the map
$\xi_M$ of \eqref{def.xi} is a p.isomorphism is analogous.  Finally,
that  $^{\vee}M \cong_p M^{\vee}$ is a formal consequence of the fact
that $\phi_M$ and $\xi_M$ are both p.isomorphisms.
\end{proof}

Although it is easy to see that $^{\vee} M$ is left noetherian over $B$,
in general there is no reason that $^{\vee} M$ should be right noetherian
over $A$.  The dualizing complex provides the extra information needed
to prove such a fact.

\begin{lemma}
\label{lem3.5}
Let $A$ and $B$ be two graded noetherian rings admitting balanced
dualizing complexes $R_A$ and $R_B$ respectively. Let $M$ be a
noetherian graded $(A,B)$-bimodule.  If $N$ is a graded noetherian right
$B$-module, then $\uExt^i_B(M,N)$ is a noetherian right $A$-module for all
$i$.
\end{lemma}

\begin{proof}
First, note that by \cite[5.1 and 4.8]{VdB},
\begin{equation}
\label{relateAB}
\RHom_{A^{\circ}}(M,R_A)\cong \R\Gamma_{{\mathfrak m}_{A^{\circ}}}(M)^*
\cong \R\Gamma_{{\mathfrak m}_B}(M)^*\cong \RHom_B(M,R_B).
\end{equation}
The functor $D(-) = \RHom_B(-, R_B)$ gives a duality from $\D(\Gr B)$
to $\D(\Gr B^{\circ})$ which restricts to a duality between the
subcategories of
complexes with finitely generated cohomology groups \cite[3.4]{Ye}. Thus
$\RHom_B(M,R_B)$ is a complex
of graded left $B$-modules with noetherian cohomologies.  Applying the
same reasoning to the duality functor
$\RHom_{A^{\circ}}(-,R_A)$ and using \eqref{relateAB}, we see that
the cohomology groups of  $\RHom_B(M,R_B)$ are also noetherian right
$A$-modules.

Now since $D(-)$ is a duality we have
$$\RHom_B(M,N)\cong \RHom_{B^\circ}(D(N),D(M)).$$
As we showed above, $D(M)$ has noetherian cohomologies on the right, and
so
$\Ext^i_{B^{\circ}}(D(N),D(M))$ also has noetherian right cohomologies.
Thus
$$\Ext^i_B(M,N) \cong \Ext^i_{B^{\circ}}(D(N),D(M))$$
is noetherian as a right $A$-module.
\end{proof}

\begin{proof}[Proof of Proposition \ref{prop3.2}]
By Lemma \ref{lem3.4} it suffices to show that $M^\vee$ and $^{\vee}M$ are
noetherian bimodules.  We only show this for $^{\vee} M$; the proof for
$M^{\vee}$ is symmetric.  It is clear that $^{\vee} M = \Hom_B(M_B, B_B)$
is a noetherian left $B$-module. By Lemma \ref{lem3.5}, $\Hom_B(M,B)$ is
also
noetherian as a right $A$-module.
\end{proof}

Let $A$ be a graded ring with balanced dualizing complex $R$, and suppose
that
$X = \Proj A$ is classically Cohen-Macaulay.  Set $M=\HS^{-(d+1)}(R)/\tau$
where $d = \cd(\Proj A)$ and $\tau$ is the torsion submodule of
$\HS^{-(d+1)}(R)$.  When we say that
the dualizing sheaf $\omega = \pi(M)$ is {\it invertible}, we mean that
$M$ is an p.invertible $(A ,A)$-bimodule in the sense defined earlier.
Now we are ready to prove Theorem~\ref{thm0.5}(b).
\begin{theorem}
\label{thm3.6} Let $A$ be a connected graded, noetherian, projectively
simple ring with a balanced dualizing complex $R$. Let $\omega$ be the
dualizing bimodule of $\Proj A$. Then $\omega$ is invertible and
${\mathcal A} = \pi(A)$ has finite injective dimension in the category
$\QGr A$.
\end{theorem}
\begin{proof} Let $M=\HS^{-(d+1)}(R)/\tau$ as in the comments before the
proof.  Then $M$ is torsionfree on both sides, so Goldie torsionfree on
both sides by Lemma~\ref{lem1.7}(a).

We need to show that $M$ is p.invertible. We already know that $M$ is
Goldie torsionfree on both sides.  Let $Q$ be the graded fraction ring of
$A$. Then $M$ is evenly localizable to $Q$ by Lemma \ref{xxlem3.1}. Since
we have already proven that $\Proj A$ is classically Cohen-Macaulay in
Proposition~\ref{prop2.1}, we know that for every $i\neq -d-1$, $\HS^i(R)$
is finite dimensional. By \cite[6.2(1)]{YZ3}, the complex $Q \otimes _A R
\otimes_A Q$ is a graded dualizing complex for the ring $Q$, and since $M$
is evenly
localizable this complex is just a shift of $M \otimes_A Q$. Since $Q$ is
graded semisimple
artinian with graded global dimension zero, a dualizing bimodule for $Q$
is a progenerator on both sides. Thus $M\otimes_A Q$ is invertible over
$Q$ and hence by definition $M$ is generically invertible. Then by
Proposition \ref{prop3.2}, $M$ is p.invertible.

By the definition of a dualizing complex $R$ has finite injective
dimension. This implies that $\pi(R)$ and hence $\omega$ has finite
injective dimension in $\QGr A$. Since $M$ is a p.invertible bimodule,
the functor $- \otimes_A M^\vee$ induces an auto-equivalence of $\QGr A$
and maps $\omega$ to $\mathcal{A}$. Therefore $\mathcal{A}$ has finite
injective dimension in $\QGr A$.
\end{proof}

\section*{Acknowledgments}
We thank Mike Artin, Lawrence Ein, Dennis Keeler, 
J\'{a}nos Koll\'{a}r, S\'{a}ndor Kov\'{a}cs, Lance Small, 
Paul Smith and Efim Zelmanov for useful discussions and 
valuable suggestions.

\providecommand{\bysame}{\leavevmode\hbox to3em{\hrulefill}\thinspace}
\providecommand{\MR}{\relax\ifhmode\unskip\space\fi MR }
\providecommand{\MRhref}[2]{%
  \href{http://www.ams.org/mathscinet-getitem?mr=#1}{#2} }
\providecommand{\href}[2]{#2}

\end{document}